\documentclass[12pt]{article}
\usepackage{amscd,amssymb,epic,amsmath}
\DeclareMathOperator{\ad}{ad}

\title{\huge\bf Degenerate principal series of quantum Harish-Chandra
modules}\author{Olga Bershtein}\date{}

\begin{document}
\newtheorem{proposition}{Proposition}
\newtheorem{lemma}{Lemma}
\newtheorem{corollary}{Corollary}

\maketitle

\centerline{\small{Institute for Low Temperature Physics \& Engineering}}
\centerline {\small{47 Lenin Ave., 61103 Kharkov, Ukraine}} \centerline
{e-mail: bershtein@ilt.kharkov.ua}

\centerline{(Received )}
\begin{abstract}

\baselineskip 2.0 pc In this paper we study a quantum analogue of a
degenerate principal series of $U_q\mathfrak{su}_{n,n}$-modules ($0<q<1$)
related to the Shilov boundary of the quantum $n \times n$-matrix unit
ball. We give necessary and sufficient conditions for the modules to be
simple and unitarizable and investigate their equivalence.

These results are q-analogues of known classical results on reducibility
and unitarizability of $SU(n,n)$-modules obtained by Johnson, Sahi, Zhang,
Howe and Tan.
\end{abstract}

\newpage

\section{Introduction}

In this paper we investigate a quantum analogue of the degenerate principal
series of representations of the algebra $U_q\mathfrak{su}_{n,n}$ related
to the Shilov boundary of the quantum $n \times n$-matrix unit ball. We
give necessary and sufficient conditions for the representations to be
irreducible and unitary.

In this work we provide q-analogues of classical results obtained by
Kenneth D.Johnson, Siddhartha Sahi, Genkai Zhang, Roger E.Howe and Eng-Chye
Tan \cite{Howe,Jons1,Jons2,Sahi,Zhang}. Another degenerate principal series
is considered in the A.Klimyk and S.Pakuliak paper \cite{KL_Pak}.

We use Bargman's approach for investigating representations (see
\cite{Barg}, where unitary strongly continuous irreducible representations
of the group $SU(1,1)$ were described). Explicit formulas for operators of
$\mathfrak{su}_{1,1}$-representations were found in a weight vectors basis
in \cite{Barg}. Results on irreducibility and unitarizability can be
obtained from the formulas as corollaries.

In the general case one need much more efforts to obtain similar formulas.
Important results in this direction were obtained by Roger Howe in
\cite{Howe}. He received certain results on irreducibility and
unitarizability of modules of the simplest degenerate principal series for
$U(m,n)$ and some other classical groups.

The Lee Soo Teck paper \cite{Lee} directly continues this Howe work. In
\cite{Lee} the degenerate principal series for $U(n,n)$ related to the
Shilov boundary of the $n \times n$-quantum ball is investigated and
answers to the same questions are obtained.

This work generalizes results from \cite{Lee} to the quantum case with
$0<q<1$. Passing to the limit as $q\mapsto 1$ one can get up to notation
the results of the above-mentioned paper.

This paper is organized as follows. In Section \ref{probl} we define the
representations $\pi_{\alpha,\beta}$ of the degenerate principal series
(see (\ref{rep})). In Section \ref{equi} we investigate the equivalence of
$\pi_{\alpha,\beta}$ (see Proposition \ref{clas_eqv}). In Section
\ref{tech} we discuss some auxiliary results $\pi_{\alpha,\beta}$. These
results will be used in the sequel. In Section \ref{reduce} we give
necessary and sufficient for $\pi_{\alpha,\beta}$ to be irreducible (see
Proposition \ref{pr_rez}). For the case $\pi_{\alpha,\beta}$ is reducible,
we describe all its irreducible subrepresentations. In Section \ref{eqv} we
find explicit formulas for intertwining operators between
$\pi_{\alpha,\beta}$ and $\pi_{\alpha,\beta}$ (see (\ref{inter_op})). In
Section \ref{unit} we investigate unitarizability of irreducible
representations of the degenerate principal series. Most of the technical
details of the proofs are contained in Appendix.

\section{Definition of the degenerate principal series of representations}
\label{probl}

Recall some concepts on geometric realizations for certain series of
representations of real semisimple Lie groups and Lie algebras.

Consider the affine algebraic group $G=SL_{2n}(\mathbb C)$ and its maximal
parabolic subgroup
\begin{equation*}
P =\left\{\left. \begin{pmatrix} A & B \\ 0 & D \end{pmatrix}|\, A,B,D \in
Mat_{n,n}(\mathbb C),\, (\det A) (\det D)=1 \right\}\right..
\end{equation*}
Then the projective variety $G \diagup P$ is isomorphic to the space
$Gr_n(\mathbb C^{2n})$ of $n$-dimensional subspaces in $\mathbb C^{2n}$.
The subgroup $K=S(GL_n(\mathbb C) \times GL_n(\mathbb C))$ acts naturally
on $G \diagup P$. Denote by $\Omega$ the open $K$-orbit. It can be easily
proved that
\begin{equation*}
 \Omega=\{L \in Gr_n(\mathbb C^{2n})
| \dim L \cap (\mathbb C^{n})_1 = \dim L \cap (\mathbb C^{n})_2 = 0\},
\end{equation*}
where $(\mathbb C^{n})_1$ and $(\mathbb C^{n})_2$ are the subspaces
generated by the elements $\{\varepsilon_1,\ldots,\varepsilon_n\}$,
$\{\varepsilon_{n+1},\ldots,\varepsilon_{2n}\}$ of the standard basis for
$\mathbb C^{2n}$, respectively. It can be verified that $\Omega$ is an
affine variety.

Set
\begin{equation*}
\mathbf t= \begin{pmatrix} t_{11} & t_{12} & \ldots & t_{1\,2n} \\ \ldots &
\ldots & \ldots & \ldots
\\t_{n\,1} & t_{n\,2} & \ldots & t_{n\,2n}
\end{pmatrix}, \quad \operatorname{rk} \,\mathbf t= n,
\end{equation*}
and $\mathbb C[\mathrm{Mat}_{n,2n}]\stackrel{\rm{def}} {=}\mathbb
C[t_{11},\ldots,t_{n\,2n}]$. Define
$$
t_{J}^{\wedge n} \stackrel{\mathrm{def}} {=}\sum_{s \in
S_n}(-1)^{l(s)}t_{1j_{s(1)}}\cdot t_{2j_{s(2)}}\cdots t_{nj_{s(n)}},
$$
where $l(s)$ is the length of permutation $s$, $J=\{j_1,\ldots,j_n\},\;1
\le j_1< \ldots <j_n \le 2n,$ and $t_{ij}$ are the matrix entries of
$\mathbf t$. The elements $t_{J}^{\wedge n}$ are called Plucker projective
"coordinates" on $Gr_n(\mathbb C^{2n})$. Denote by $\mathbb C[\mathrm
{Pl}_{n,2n}] \subset \mathbb C[\mathrm{Mat}_{n,2n}]$ the subalgebra
generated by all $t^{\wedge n}_J$.

Consider the algebra $\mathbb C[\Omega]$ of regular functions on $\Omega$.
Let us introduce some notation. Set $t\stackrel{\rm def}=t^{\wedge
n}_{\{n+1,\ldots,2n\}}$ and
\begin{align*}
z_a^b=t^{-1}t^{\wedge n}_{J_{a\,b}}, \quad &a,b=1,\ldots,n, \quad
\text{where} \quad J_{a\,b}=\{n+1,\ldots,2n\} \backslash \{2n+1-b\} \cup
\{a\};
\\& \mathbf z= \begin{pmatrix} z_1^1 & \ldots & z_1^n \\ \ldots & \ldots & \ldots
\\ z_n^1 & \ldots &
z_n^n \end{pmatrix}, \quad \det \mathbf z =\det
\begin{pmatrix} z_1^1 & \ldots & z_1^n
\\ \ldots & \ldots & \ldots
\\ z_n^1 & \ldots & z_n^n \end{pmatrix}.
\end{align*}
Then the algebra $\mathbb C[\Omega]$ is canonically isomorphic to the
localization of the algebra $\mathbb C[\mathrm{Mat}_n] \stackrel {\rm
def}{=} \mathbb C[z_1^1,\ldots,z_n^n]$ with respect to the multiplicative
set $(\det \mathbf z)^{\mathbb Z_+}$. The vector space $\mathbb
C[\Omega]=\mathbb C[\mathrm{Mat}_n]_{\det \mathbf z}$ can be naturally
equipped with an $\mathfrak{sl}_{2n}$-module structure and a $K$-module
structure, and these structures are compatible (see \cite{Wall}).

Therefore the action of the universal enveloping algebra
$U\mathfrak{sl}_{2n}$ in the vector space $\mathbb
C[\mathrm{Mat}_n]_{\det\mathbf z}$ is well defined. Moreover, the
$U\mathfrak{sl}_{2n}$-action in the localization of the algebra $\mathbb
C[\mathrm{Pl}_{n,2n}]$ with respect to the multiplicative set $t^{\mathbb
Z_+}$ is well defined. Hence the $U\mathfrak{sl}_{2n}$-action in the
space\footnote{They are spaces of sections of homogeneous vector bundles
over $\Omega$. We pass from $\alpha,\beta \in \mathbb{Z}$ to $\alpha,\beta
\in \mathbb{R}$ standardly.} $\mathbb C[\mathrm{Mat}_n]_{\det\mathbf z}
\cdot (\det\mathbf z)^{\alpha} t^{\beta}$ is well defined for each
$\alpha,\beta \in \mathbb{Z}$.

\medskip

Now let us pass to the quantum case. Everywhere in the sequel $q \in
(0,1)$, $\mathbb C$ is the ground field and all algebras are unital.

Denote by $U_q\mathfrak{sl}_{2n}$ the algebra defined by its generators
$\{E_i,\:F_i,\:K_i,\:K_i^{-1}\}_{i=1}^{2n-1}$ and the relations
$$
K_iK_j \,=\,K_jK_i,\quad K_iK_i^{-1}\,=\,K_i^{-1}K_i \,=\,1;
$$
$$
K_iE_i \,=\,q^{2}E_iK_i,\quad K_iF_i \,=\,q^{-2}F_iK_i;
$$
$$
E_iF_j \,-\,F_jE_i \:=\:\delta_{ij}(K_i-K_i^{-1})/(q-q^{-1});
$$
$$
K_iE_j \,=\,q^{-1}E_jK_i,\, K_iF_j \,=\,qF_jK_i, \quad |i-j|=1;
$$
$$
E_i^2E_j\,-\,(q+q^{-1})E_iE_jE_i \,+\,E_jE_i^2 \:=\:0,\quad |i-j| \,=\,1;
$$
$$
F_i^2F_j \,-\,(q+q^{-1})F_iF_jF_i \,+\,F_jF_i^2\:=\:0,\quad |i-j|\,=\,1;
$$
$$
K_iE_j \,-\,E_jK_i\,=\, K_iF_j \,-\,F_jK_i,\,=\,
E_iE_j\,-\,E_jE_i\,=\,F_iF_j\,-\,F_jF_i\,=\,0, \quad |i-j|>1.
$$

We equip $U_q\mathfrak{sl}_{2n}$ with the standard Hopf algebra structure.
The comultiplication, the counit and the antipode are defined by their
actions on the generators:
\begin{align*}
\triangle{E_j}&=E_j \otimes 1+K_j \otimes E_j, & \varepsilon(E_j)&=0,
& S(E_j)=&-K_j^{-1}E_j,\\
\triangle{F_j}&=F_j \otimes K_j^{-1}+1 \otimes F_j, & \varepsilon(F_j)&=0,
& S(F_j)=&-F_jK_j,\\
\triangle{K_j}&=K_j \otimes K_j, & \varepsilon(K_j)&=1, & S(K_j)=&K_j^{-1}
\end{align*}
for all $j=1,\ldots,2n-1$.

The algebra $\mathbb{C}[\mathrm{Mat}_{n,2n}]_q$ of polynomials on the
quantum $n\times2n$-matrix space is defined by its generators
$\{t_{ij}\}_{i=1,\ldots,n;j=1,\ldots,2n}$ and the relations (cf.
\cite{Drin_con})
\begin{align}\label{q_rel}
&t_{ik}t_{jk}=qt_{jk}t_{ik},\quad t_{ki}t_{kj}=qt_{kj}t_{ki},&i<j, \notag
\\ &t_{ij}t_{kl}=t_{kl}t_{ij},&i<k \; \& \; j>l,
\\ &t_{ij}t_{kl}-t_{kl}t_{ij}=(q-q^{-1})t_{ik}t_{jl},&i<k \; \& \; j<l.
\notag
\end{align}
Define q-minors as follows:
\begin{equation}\label{q-min}
t_{IJ}^{\wedge k}\stackrel{\mathrm{def}}{=}\sum_{s \in
S_k}(-q)^{l(s)}t_{i_1j_{s(1)}}\cdots t_{i_kj_{s(k)}},
\end{equation}
for any $I=\{i_1,\ldots,i_k\},\;1 \le i_1<\ldots<i_k \le n,$
$J=\{j_1,\ldots,j_k\}, \;1 \le j_1<\ldots<j_k \le 2n$; here $l(s)$ denotes
the length of permutation $s$.

Consider the algebra $\mathbb{C}[\mathrm{Pl}_{n,2n}]_{q} \subset
\mathbb{C}[\mathrm{Mat}_{n,2n}]_{q}$ generated by all q-minors $t^{\wedge
n}_{\{1,\ldots,n\}J}$, $\operatorname{card} J=n$. It is equipped with the
standard $U_q \mathfrak{sl}_n^{op} \otimes U_q \mathfrak{sl}_{2n}$-module
algebra structure.\footnote{$U_q \mathfrak{sl}_n^{op}$ is a Hopf algebra
with the same multiplication and the opposite comultiplication.} It is easy
to show that the $U_q \mathfrak{sl}_n^{op}$-structure can be reconstructed
from the below equalities:
\begin{align*}
K_l \,t_{ij}=&\begin{cases} q^{-1}t_{ij}, & l=i,
\\qt_{ij}, & l=i-1,
\\ 0, &\text{otherwise};
\end{cases}
\\ E_l \,t_{ij}=q^{-1/2} \cdot \begin{cases}  t_{(i+1)j}, & l=i,
\\0, &\text{otherwise};
\end{cases}
& \quad F_l \,t_{ij}=q^{1/2} \cdot \begin{cases} t_{(i-1)j}, & l=i-1,
\\0, &\text{otherwise}.
\end{cases}
\end{align*}

The element $t\stackrel{\rm def}{=} t^{\wedge
n}_{\{1,2,\ldots,n\}\{n+1,n+2,\ldots,2n\}}$ quasi-commutes with $t_{ij}$
for all $i=1,\ldots,n$, $j=1,\ldots,2n$ and is
$U_q\mathfrak{sl}^{op}_n$-invariant.

Denote by $\mathbb{C}[\mathrm{Pl}_{n,2n}]_{q,t}$ the localization of the
algebra $\mathbb{C}[\mathrm{Pl}_{n,2n}]_q$ with respect to the
multiplicative system $t^{\mathbb Z_+}$. Introduce q-analogues of
coordinates on $\Omega$ as follows:
\begin{equation}\label{emb}
z_a^b \stackrel{\mathrm{def}}{=} t^{-1}t^{\wedge n}_{\{1,2,\ldots,n \}J_{a
\,b}},
\end{equation}
where $J_{a \,b}=\{n+1,n+2,\ldots,2n \}\setminus \{2n+1-b \}\cup \{a \}.$

The defining relations for the subalgebra generated by the elements $z_a^b$
are obtained in \cite{V_qmbdic}:
\begin{align*}
&z_{a}^{b_1}z_{a}^{b_2}=qz_{a}^{b_2}z_{a}^{b_1}, &b_1<b_2,
\\&z_{a_1}^{b}z_{a_2}^{b}=qz_{a_2}^{b}z_{a_1}^{b}, &a_1<a_2,
\\&z_{a_1}^{b_1}z_{a_2}^{b_2}=z_{a_2}^{b_2}z_{a_1}^{b_1},
&b_1<b_2 \, \&\, a_1>a_2,
\\&z_{a_1}^{b_1}z_{a_2}^{b_2}-z_{a_2}^{b_2}z_{a_1}^{b_1}=
(q-q^{-1})z_{a_1}^{b_2}z_{a_2}^{b_1}, &b_1<b_2 \, \&\, a_1<a_2.
\end{align*}
(For the special case $n=2$ see the Noumi paper \cite{Noumi}.)

It can be checked easily that $zt=qtz$ for any $z \in \{z_a^b|
a,\,b=1,\ldots,n\}$.

%\begin{proposition}
%There is a unique extension of $U_q \mathfrak{sl}_{2n}$-module algebra
%structure from $\mathbb{C}[\mathrm{Pl}_{n,2n}]_q$ onto
%$\mathbb{C}[\mathrm{Pl}_{n,2n}]_{q,t}$.
%\end{proposition}
%{\bf Proof.} The uniqueness is obvious. Then for each $f \in
%\mathbb{C}[\mathrm{Pl}_{n,2n}]_{q,t}$ there exists an integer $k \in
%\mathbb N$ such that $ft^k \in \mathbb{C}[\mathrm{Pl}_{n,2n}]_q$. This
%allows one to extend the $U_q \mathfrak{sl}_{2n}$-module algebra structure
%from $\mathbb{C}[\mathrm{Pl}_{n,2n}]_q$ onto
%$\mathbb{C}[\mathrm{Pl}_{n,2n}]_{q,t}$. We just have to verify that this
%extension is well defined. For that reason we mention that for any $\xi \in
%U_q \mathfrak{sl}_{2n}$, $f \in \mathbb{C}[\mathrm{Pl}_{n,2n}]_{q,t}$ there
%is a unique Laurent polynomial $p_{f,\xi}$ of the variable $u=q^{k}$ with
%coefficients from $\mathbb{C}[\mathrm{Pl}_{n,2n}]_{q,t}$ such that
%$p_{f,\xi}(q^k)t^k=\xi\cdot(f t^k)$, and such polynomials satisfy the
%conditions $p_{f,\xi}(q^l)t=p_{ft,\xi}(q^{l-1})$ for $l \in \mathbb N$.
%\hfill $\square$

It can be proved that for any $\xi \in U_q \mathfrak{sl}_{2n}$, $f \in
\mathbb{C}[\mathrm{Pl}_{n,2n}]_{q,t}$ there is a unique Laurent polynomial
$p_{f,\xi}$ of the variable $u=q^{k}$ with coefficients in
$\mathbb{C}[\mathrm{Pl}_{n,2n}]_{q,t}$ such that $p_{f,\xi}(q^k)=\xi\cdot(f
t^k)t^{-k}$. This allows one to prove the existence of an extension of $U_q
\mathfrak{sl}_{2n}$-module algebra structure onto
$\mathbb{C}[\mathrm{Pl}_{n,2n}]_{q,t}$ (see \cite{Geom_real}).

The subalgebra generated by $z_a^b$ is the algebra $\mathbb
C[\mathrm{Mat}_n]_q$ of ''polynomials on the quantum $n \times n$-matrix
space '' (cf. (\ref{q_rel})). The algebra $\mathbb C[\mathrm{Mat}_n]_q$ is
a $U_q\mathfrak{sl}_{2n}$-module subalgebra of the
$U_q\mathfrak{sl}_{2n}$-module algebra
$\mathbb{C}[\mathrm{Pl}_{n,2n}]_{q,t}$ (see \cite{V_ftqmbii}).
\begin{proposition}\label{uqslmat}
(\cite{V_ftqmbii}) For all ${a,b=1,\ldots,n}$
\begin{align*}
K_n^{\pm 1}z_a^b=&
\begin{cases}
q^{\pm 2}z_a^b,&a=n \;\&\;b=n
\\ q^{\pm 1}z_a^b,&a=n \;\&\;b \ne n \;\text{or}\;
a \ne n \;\&\; b=n
\\ z_a^b,&\text{\rm otherwise},
\end{cases}
\\ F_nz_a^b=q^{1/2}\cdot
\begin{cases}
1,& a=n \;\& \;b=n
\\ 0,&\text{\rm otherwise},
\end{cases}\quad
& \quad E_nz_a^b=-q^{1/2}\cdot
\begin{cases}
q^{-1}z_a^nz_n^b,&a \ne n \;\&\;b \ne n
\\ (z_n^n)^2,& a=n \;\&\;b=n
\\ z_n^nz_a^{b},&\text{\rm otherwise}
\end{cases}
\end{align*}
and for all $k \ne n$ we have
\begin{align*}
\!\!\!\!\!\!K_k^{\pm 1}z_a^b=
\begin{cases}
\!q^{\pm 1}z_a^b,\; k<n \;\&\;a=k \;\text{or} \;k>n \;\&\;b=2n-k,
\\ \!q^{\mp 1}z_a^b,\; k<n \;\&\;a=k+1 \;\text{or} \;k>n \;\&\;b=2n-k+1,
\\ \!z_a^b,\; \text{\rm otherwise},
\end{cases}
\\ \!\!\!F_kz_a^b=\!q^{1/2}\!\cdot
\begin{cases}
\!z_{a+1}^b, \;k<n \;\&\;a=k,
\\ \!z_a^{b+1}, \;k>n \;\&\;b=2n-k,
\\ \!0, \;\; \text{\rm otherwise},
\end{cases}
\!\!\! E_kz_a^b=\!q^{-1/2}\!\cdot
\begin{cases}
\!z_{a-1}^b, \;k<n \;\&\; a=k+1,
\\ \!z_a^{b-1}, \;k>n \;\&\;b=2n-k+1,
\\ \!0, \;\;\;\;\text{\rm otherwise}.\hfill \square
\end{cases}
\end{align*}
\end{proposition}

In the sequel we use the following notation for q-minors
\begin{equation}\label{q_minz}
\mathbf z ^{\wedge k \{b_1,\ldots,b_k\}}_{\,\,\,\,\,\,\,\{a_1,\ldots,a_k\}}
\stackrel {\rm def}{=} \sum_{s \in S_k} (-q)^{l(s)}
z^{b_{s(1)}}_{a_1}\ldots z^{b_{s(k)}}_{a_k},
\end{equation}
where $a_1 < \ldots < a_k,\,\, b_1<\ldots<b_k$. It is known that the
element $\det_q \mathbf z\stackrel{\mathrm{def}}{=}\mathbf {z}
_{\,\,\,\,\,\,\,\{1,\ldots,n\}}^{\wedge n\{1,\ldots,n\}}$ belongs to the
center of $\mathbb C[\mathrm{Mat}_n]_q$ and $\mathbb C[\mathrm{Mat}_n]_q$
has no zero divisors.

Denote by $\mathbb C[\mathrm{Mat}_n]_{q,\det_q \mathbf z}$ the localization
of the algebra $\mathbb C[\mathrm{Mat}_n]_q$ with respect to the
multiplicative system $(\det_q \mathbf z)^{\mathbb Z_+}$. We consider
$\mathbb C[\mathrm{Mat}_n]_{q,\det_q \mathbf z}$ as a q-analogue of the
space of regular functions on the open orbit $\Omega$. Let
$\widetilde{t}=t^{\wedge n}_{\{1,\ldots,n\} \{1,\ldots,n\}}$. Since $\det_q
\mathbf z = t^{-1}\widetilde{t}$, we see that the algebra $\mathbb
C[\mathrm{Mat}_n]_{q,\det_q \mathbf z}$ is a $U_q\mathfrak{sl}_{2n}$-module
subalgebra of the $U_q\mathfrak{sl}_{2n}$-module algebra
$\mathbb{C}[\mathrm{Pl}_{n,2n}]_{q,t \cdot \widetilde{t}}$. (As above, to
verify that the extension is well defined we use the following fact: for
all $\xi \in U_q\mathfrak{sl}_{2n}, \, f \in V$ the vector valued function
$\xi \cdot (f (\det_q \mathbf z)^{k}) (\det_q \mathbf z)^{-k}$ is a Laurent
polynomial of the variable $u=q^{k}$.)

Denote by $V$ the vector space $\mathbb C[\mathrm{Mat}_n]_{q,\det_q \mathbf
z}$. Assume first that $\alpha,\beta \in \mathbb Z$. Define a
representation $\pi_{\alpha,\beta}:U_q\mathfrak{sl}_{2n} \rightarrow {\rm
End} V$ as follows:
\begin{equation}\label{rep}
\pi_{\alpha,\beta}(\xi)f=(\xi \cdot (f (\widetilde{t})^{\alpha} t^{\beta}))
t^{-\beta} (\widetilde{t})^{-\alpha} = (\xi \cdot (f ({\rm det}_q \mathbf
z)^{\alpha} t^{\beta+\alpha})) t^{-\alpha-\beta} ({\rm det}_q \mathbf
z)^{-\alpha}
\end{equation}
for every $\xi \in U_q\mathfrak{sl}_{2n},f \in V$. For each $\lambda \in
\mathbb Z$ we have
$$
E_j t^{\lambda}=0,\,\, F_jt^{\lambda}=0,\,\, K_jt^{\lambda}=t^{\lambda},
\quad j=1,\ldots,2n-1, \; j\ne n
$$
$$
E_nt^\lambda=q^{-3/2}\frac{1-q^{-2\lambda}}{1-q^{-2}}z^n_n t^\lambda,
\quad F_nt^\lambda=0, \quad K_n^{\pm 1}t^\lambda=q^{\mp
\lambda}t^\lambda,
$$
$$
E_j ({\rm det}_q \mathbf z)^\lambda=0,\,\, F_j ({\rm det}_q \mathbf
z)^\lambda=0,\,\, K_j ({\rm det}_q \mathbf z)^\lambda=({\rm det}_q \mathbf
z)^\lambda, \quad j=1,\ldots,2n-1, \; j\ne n
$$
$$
K_n^{\pm 1}(({\rm det}_q \mathbf z)^\lambda)=q^{\pm 2 \lambda}({\rm det}_q
\mathbf z)^\lambda, \quad E_n(({\rm det}_q \mathbf z)^\lambda) = -q^{1/2}
\frac{1-q^{2\lambda}}{1-q^2} z^n_n ({\rm det}_q \mathbf z)^\lambda,
$$
$$
F_n(({\rm det}_q \mathbf z)^\lambda) =
q^{1/2}\frac{1-q^{-2\lambda}}{1-q^{-2}} z^{\wedge
n-1}_{\{1,\ldots,n-1\}\{1,\ldots,n-1\}} ({\rm det}_q \mathbf
z)^{\lambda-1},\quad \lambda \ne 0.
$$

From these equalities we see that for each $\xi \in U_q\mathfrak{sl}_{2n}$,
$f \in V$ the vector valued function $p_{f,\xi}(q^{\alpha},q^{\beta})
\stackrel{\rm def}=\pi_{\alpha,\beta}(\xi)(f)$ is a Laurent polynomial of
the variables $q^{\alpha},\,q^{\beta}$. These Laurent polynomials are
defined by their values on the set $\{(q^{\alpha},q^{\beta})|
\,\alpha,\,\beta \in \mathbb{Z}\}$ and deliver the canonical "analytic
continuation" for $\pi_{\alpha,\beta}(\xi)(f)$ to $(\alpha,\beta) \in
\mathbb{C}^{2}$.

Let $(\alpha,\beta) \in \mathbb{C}^{2}$. Define a representation
$\pi_{\alpha,\beta}(\xi)(f) \stackrel{\rm
def}=p_{f,\xi}(q^{\alpha},q^{\beta})$. Indeed, to prove that the
representation $\pi_{\alpha,\beta}$ is well defined for
$(q^{\alpha},\,q^{\beta}) \in \mathbb{C}^{2}$ it is sufficient to verify
some identities for Laurent polynomials. These identities are correct for
$\alpha,\,\beta \in \mathbb{Z}$.

Introduce a ''deformation parameter'' $h$ by the equality $q=e^{-h/2}$.
Clearly, if $\alpha_1=\alpha_2+i\frac{2\pi}{h}$ and
$\beta_1=\beta_2+i\frac{2\pi}{h}$, then
$\pi_{\alpha_1,\beta_1}=\pi_{\alpha_2,\beta_2}$. Then it is enough to
consider $\alpha,\,\beta \in D$, where
\begin{equation*}
D=\{ \alpha,\beta \in \mathbb C \, |\, 0 \leq \operatorname{Im} \alpha <
\frac{2\pi}{h}, \, 0 \leq \operatorname{Im} \beta < \frac{2\pi}{h}\}.
\end{equation*}

Recall that a representation $\rho:U_q \mathfrak{sl}_{2n} \rightarrow {\rm
End} W$ is called {\it weight} if the representation space $W$ decomposes
as follows:
\begin{align*}
&W=\bigoplus_{\boldsymbol \lambda} W_{\boldsymbol \lambda}, \quad
\text{where} \quad \boldsymbol \lambda =(\lambda_1,\ldots,\lambda_{2n-1})
\in \mathbb Z^{2n-1},
\\& W_{\boldsymbol \lambda}=\{v \in W |\,\rho
(K_j^{\pm})v = q^{\pm \lambda_j}v,j=1,\ldots,2n-1 \}.
\end{align*}
The subspace $W_\lambda$ is called weight subspace with weight $\lambda$.
In the sequel we will consider only weight representations. It is clear
that $\pi_{\alpha,\beta}$ is a weight representation if and only if
$q^{\alpha-\beta} \in q^{\mathbb Z}$.

Let $W$ be a weight $U_q\mathfrak{sl}_{2n}$-module. Define operators $H_i$
for $i=1,\ldots,2n-1$ by the formula $H_i|_{W_{\lambda}}=\lambda_i$.

\section {Equivalence of the representations}\label{equi}

Recall that $q=e^{-h/2}$. For any complex $\alpha,\beta$ such that $0 \leq
\operatorname{Im} \alpha < \frac{2\pi}{h}, \, 0 \leq \operatorname{Im}
\beta < \frac{2\pi}{h}$, the statements $\alpha-\beta \in \mathbb Z$ and
$q^{\alpha-\beta} \in q^{\mathbb Z}$ are equivalent.
\begin{proposition} \label{pr_eqv}
If $\alpha,\beta \not \in \mathbb Z$, then the representations
$\pi_{\alpha,\beta}$
 and $\pi_{-n-\beta,-n-\alpha}$ are equivalent.
\end{proposition}
The proof is reduced to explicit formulas for the intertwining operators.
It is given in Section \ref{eqv}.

If $\alpha,\beta \in \mathbb Z$, then the representations
$\pi_{\alpha,\beta}$ and $\pi_{-n-\beta,-n-\alpha}$ are not equivalent.
This fact follows from the statement that only one of the representations
$\pi_{\alpha,\beta}$ and $\pi_{-n-\beta,-n-\alpha}$ for integral
$\alpha,\beta$ has a finite dimensional subrepresentation. An explanation
of this fact is given in the end of Section \ref{reduce}.

The representations $\pi_{\alpha,\beta}$ and $\pi_{\alpha-1,\beta+1}$ are
equivalent for all $\alpha,\beta$. The corresponding intertwining operator
$T: V \rightarrow V$ is defined as follows: for every $f \in V=\mathbb
C[\mathrm{Mat}_n]_{q,\det_q \mathbf z}$ $T(f)=f(\det_q \mathbf z)^{-1}.$
Indeed, since for each $f \in V,\,\xi \in U_q \mathfrak{sl}_{2n}$
\begin{multline*}
\pi_{\alpha-1,\beta+1}(\xi)(f)=(\xi \cdot (f ({\rm det}_q \mathbf
z)^{\alpha-1} t^{\beta+\alpha})) t^{-\alpha-\beta} ({\rm det}_q \mathbf
z)^{1-\alpha}=
\\ (\xi \cdot (f ({\rm det}_q \mathbf z)^{-1} ({\rm det}_q \mathbf z)^{\alpha}
t^{\beta+\alpha})) t^{-\alpha-\beta} ({\rm det}_q \mathbf z)^{-\alpha}
({\rm det}_q \mathbf z)= \pi_{\alpha,\beta}(\xi)(f({\rm det}_q \mathbf
z)^{-1}){\rm det}_q \mathbf z,
\end{multline*}
we see that $T$ intertwines the representations $\pi_{\alpha,\beta}$ and
$\pi_{\alpha-1,\beta+1}$. Therefore without loss of generality we can
assume that $\alpha,\beta \in \mathcal{D}$, where
\begin{equation}\label{lim} \mathcal{D}=\{ (\alpha,\beta) \in \mathbb
C \, |\, \alpha-\beta \in \{0,1\},\, 0 \leq \operatorname{Im} \alpha <
\frac{2\pi}{h}, 0 \leq \operatorname{Im} \beta < \frac{2\pi}{h}\}.
\end{equation}
Let us introduce an equivalence relation on $\mathcal{D}$. The equivalence
class of $(\alpha,\beta)$ consists of one point for $\alpha,\beta \in
\mathbb Z$ and from two points for $\alpha,\beta \not \in \mathbb Z$:
\begin{equation*}
(\alpha_1,\beta_1) \sim (\alpha_2,\beta_2), \quad \text{iff} \quad
\begin{cases} \alpha_1=-n-\beta_2,\,\beta_1=-n-\alpha_2 \,\,\text{for} \,\,
\operatorname{Im} \alpha_1 = \operatorname{Im} \alpha_2 =0,
\\ \alpha_1=\frac{2 \pi i}{h}-n-\beta_2,\,\beta_1=\frac{2 \pi i}{h}-n-\alpha_2,
\,\, \text{otherwise}.
\end{cases}
\end{equation*}
\begin{proposition}\label{clas_eqv}
The set of equivalence classes $\mathcal{D} \diagup \sim$ is in the
one-to-one correspondence $(\alpha,\beta) \mapsto \pi_{\alpha,\beta}$ with
the set of equivalence classes of the representations of the degenerate
principal series.
\end{proposition}
{\bf Proof.} By the above, each representation of the degenerate principal
series is equivalent to the representation $\pi_{\alpha,\beta}$ for some
$(\alpha,\beta) \in \mathcal{D}$.

Prove that the representations $\pi_{\alpha_1,\beta_1}$ and
$\pi_{\alpha_2,\beta_2}$, with $(\alpha_1,\beta_1),(\alpha_2,\beta_2) \in
\mathcal{D}$, are equivalent if and only if $(\alpha_1,\beta_1) \sim
(\alpha_2,\beta_2)$. For that we calculate the action of a central element
$C \in U_q\mathfrak{sl}_{2n}^{ext}$ (see \cite{Klimyk} for the definition).
It can be proved that $\pi_{\alpha,\beta}(C)$ is a scalar operator for all
$\alpha,\beta \in \mathcal{D}$.

From \cite{Drin} it follows that there exists a unique central element $C$
which acts on the $U_q\mathfrak{sl}_{2n}$-highest vector $v^{high}$ with
weight $\lambda$ %(see definition in the end of Section \ref{probl})
as follows:
\begin{equation*}
C({v^{high}})=\sum \limits_{j=0}^{2n-1} q^{-2(\mu_j,\lambda+\rho)}v^{high},
\end{equation*}
where $\mu_0=\varpi_1, \;\mu_j=-\varpi_j+\varpi_{j+1}$ for
$j=1,\ldots,2n-2$, $\mu_{2n-1}=-\varpi_{2n-1}$, $\varpi_j$ are the
fundamental weights, $2\rho$ is the sum of positive roots of the Lie
algebra $\mathfrak{sl}_{2n}$, and we choose the invariant scalar product
such that $(\alpha,\,\alpha)=2$ for any simple root $\alpha$.

First let $\alpha,\beta$ be integers. It can be proved that
\begin{equation*}
\pi_{\alpha,\beta}(C)({\rm det}_q \mathbf z)^{\beta}=
4\operatorname{ch}\frac{h}{2} (\alpha+\beta+n) (\sum_{j=0}^{n-1}
\operatorname{ch} \frac{h}{2}j) ({\rm det}_q \mathbf z)^{\beta}.
\end{equation*}
Hence $\pi_{\alpha,\beta}(C)=4\operatorname{ch}\frac{h}{2} (\alpha+\beta+n)
(\sum_{j=0}^{n-1} \operatorname{ch} \frac{h}{2}j)\cdot \mathrm{Id}$ for all
$(\alpha,\beta) \in \mathcal{D}$.

Suppose that $\pi_{\alpha_1,\beta_1}$ and $\pi_{\alpha_2,\beta_2}$ are
equivalent. Equivalent representations have the same weight lattice.
Therefore $(\alpha_1-\beta_1)-(\alpha_2-\beta_2) \in 2\mathbb Z$. Since
$(\alpha_1,\beta_1),(\alpha_2,\beta_2) \in \mathcal{D}$, we see that
$(\alpha_1-\beta_1)-(\alpha_2-\beta_2)=0$.

Then the equivalent representations $\pi_{\alpha_1,\beta_1}$ and
$\pi_{\alpha_2,\beta_2}$ have the same values of central characters, which
means that
$$
(\operatorname{ch}\frac{h}{2}(\alpha_1+\beta_1+n)-
\operatorname{ch}\frac{h}{2}(\alpha_2+\beta_2+n)) \sum_{j=0}^{n-1}
\operatorname{ch} \frac{h}{2}j=0
$$
Since $0 \leq \operatorname{Im} \alpha_1 < \frac{2\pi}{h}$, $0 \leq
\operatorname{Im} \beta_1 < \frac{2\pi}{h}$, $0 \leq \operatorname{Im}
\alpha_2 < \frac{2\pi}{h}$, $0 \leq \operatorname{Im} \beta_2 <
\frac{2\pi}{h}$, we have that $\alpha_1+\beta_1=\alpha_2+\beta_2$, or
$\alpha_1+\beta_1=-\alpha_2-\beta_2-2n$, or
$\alpha_1+\beta_1=-\alpha_2-\beta_2-2n-\frac{4\pi i}{h}$. If
$\alpha_1+\beta_1=\alpha_2+\beta_2$, then $\alpha_1=\alpha_2$ and $\beta_1=
\beta_2$. For any fixed non-integral $\alpha_1,\beta_1$ there is a unique
pair $(\alpha_2,\beta_2) \in \mathcal{D}$ such that either
$\alpha_1+\beta_1=-\alpha_2-\beta_2-2n$ or
$\alpha_1+\beta_1=-\alpha_2-\beta_2-2n-\frac{4\pi i}{h}$, and
$(\alpha_1,\beta_1) \sim (\alpha_2,\beta_2)$. Although for integral
parameters $\pi_{\alpha_1,\beta_1}$ and $\pi_{\alpha_2,\beta_2}$ are not
equivalent, because the only one of them has a finite-dimensional
subrepresentation. This can be deduced from Corollary \ref{cor_4}. Thus
each equivalence class in $\mathcal{D}$ is assigned to a unique equivalence
class of the representations of the degenerate principal series
$\pi_{\alpha,\beta}$.\hfill $\square$

\section {Auxiliary statements on $\pi_{\alpha,\beta}$-structure}
\label{tech}

In this section we describe some necessary technical results, that will be
useful in the sequel.

Everywhere in this section we assume that $n>1$. However, Propositions
\ref{K_reduce}, \ref{pr_2}, and \ref{pr_3} and Corollaries \ref{cor_2} and
\ref{cor_3} are still sensible and correct for $n=1$.

Let $U_q \mathfrak{k}_{ss}\subset U_q \mathfrak{sl}_{2n}$ be the Hopf
subalgebra generated by $E_j, F_j, K_j^{\pm 1}$, $j=1,\ldots,2n-1$, $j \neq
n$ and $U_q \mathfrak{k}\subset U_q \mathfrak{sl}_{2n}$ be the Hopf
subalgebra generated by $K_n^{\pm 1}$ and $U_q \mathfrak{k}_{ss}$.

Note that $\pi_{\alpha,\beta}|_{U_q \mathfrak{k}_{ss}}$ does not depend on
$\alpha,\beta$. The following preliminary result on reducibility of
$\pi_{\alpha,\beta}$ is well known in the classical case. For brevity,
set\footnote{Note that, obviously, $\mathbf z^{\wedge n}=\det_q \mathbf
z$.}
$$
\mathbf z^{\wedge k} = \mathbf z^{\wedge k
\{1,\ldots,k\}}_{\,\,\,\,\,\,\{1,\ldots,k\}}.
$$

Introduce the following notation: $\widehat{K}=\{\overline{\mathbf
k}=(k_1,\ldots,k_n) \in \mathbb Z^n|\,k_1 \geq k_2 \geq \ldots \geq k_n\}$,
$\mathbf{e_j}=(0,\ldots,\overset{j}{1},\ldots,0) \in \mathbb Z^n$.

\begin{proposition} \label{K_reduce}
The representation space $V$ for $\pi_{\alpha,\beta}$ splits into a sum of
simple pairwise non-isomorphic $U_q \mathfrak{k}$-modules as
follows:\footnote{These isotypic components are $U_q
\mathfrak{k}_{ss}$-isomorphic. However, they are not $U_q
\mathfrak{k}$-isomorphic, since the action of $\pi_{\alpha,\beta}(K_n)$
depends on $\alpha,\beta$.}
\begin{equation*}
V=\bigoplus_{\overline{\mathbf k} \in \widehat{K}} V_{\overline{\mathbf
k}}, \quad \text{with} \quad V_{\overline{\mathbf
k}}=\pi_{\alpha,\beta}(U_q\mathfrak{k}) \cdot v^h_{\overline{\mathbf k}}
\quad \text{and} \quad v^h_{\overline{\mathbf k}}=(z^{\wedge
1})^{k_1-k_2}\ldots(z^{\wedge n-1})^{k_{n-1}-k_n}(z^{\wedge n})^{k_n}.
\end{equation*}
\end{proposition}

{\bf Proof.} Consider the filtration $V=\bigcup_{k=0}^{\infty} V^{(k)}$
with $V^{(k)}=\mathbb C [\mathrm{Mat}_n]_q\cdot(\det_q \mathbf z)^{-k}$. It
is sufficient to prove that
\begin{equation*}
V^{(k)}=\bigoplus \limits_{k_n \geq -k} V_{\overline{\mathbf k}}.
\end{equation*}

Equip the vector space $V^{(k)}$ with the natural grading
$V^{(k)}=\bigoplus \limits_{j=-nk}^{\infty} (V^{(k)})_j$ as follows:
$(V^{(k)})_j=\{v \in V^{(k)}|\,K_0v=q^{2j}v\}$, with $K_0 \stackrel {\rm
def}{=} K_1 K_2^2K_3^3 \ldots K_n^nK_{n+1}^{n-1}\ldots K_{2n-2}^2K_{2n-1}$.
Therefore we must prove that
\begin{align}\label{l}
(\mathbb C [\mathrm{Mat}_n]_q\cdot({\rm det}_q \mathbf z)^{-k})_j =
\bigoplus \limits_{\substack{k_n \geq -k,\\
k_1+\ldots+k_n=j}} V_{\overline{\mathbf k}}.
\end{align}

For $k=0$ statement (\ref{l}) means that $(\mathbb C [\mathrm{Mat}_n]_q)_j
= \bigoplus \limits_{k_n \geq 0 ,\,k_1+\ldots+k_n=j} V_{\overline{\mathbf
k}}.$

First, the dimensions of homogeneous components $\mathbb C
[\mathrm{Mat}_n]_{q,j}$ of the standardly graded algebra $\mathbb C
[\mathrm{Mat}_n]_q$ are equal to the dimensions in the classical case:
\[
\dim \, \mathbb C [\mathrm{Mat}_n]_{q,j} \,=\,\begin{pmatrix} n^2+j-1 \\ j
\end{pmatrix}
\]
(it can be easily proved via the Bergman diamond lemma \cite{Klimyk}, sec.
4.1.5). Secondly, the dimensions of the $U_q\mathfrak{k}$-modules
$V_{\overline{\mathbf k}}$ are equal to the classical ones (this follows
from results of quantum groups theory \cite{Jantzen}, chap. 5). Thirdly,
there is the well-known Hua result on the coincidence of the dimensions
$\mathbb C [\mathrm{Mat}_n]_j$ and $\bigoplus
\limits_{\substack{k_n \geq 0, \\
k_1+k_2+\ldots+k_n=j}} V_{\overline{\mathbf k}}$ in the classical case
\cite{Hua}. Hence,
\begin{equation*}
\dim (\mathbb C [\mathrm{Mat}_n]_q)_j = \sum \limits_{\substack{k_n \geq 0
,\\k_1+\ldots+k_n=j}} \dim V_{\overline{\mathbf k}},
\end{equation*}
and, finally,
\begin{equation*}
(\mathbb C [\mathrm{Mat}_n]_q)_j = \bigoplus
\limits_{\substack{k_n \geq 0 ,\\
k_1+\ldots+k_n=j}} V_{\overline{\mathbf k}}.
\end{equation*}

For $k>0$ one has
\begin{multline*}
(\mathbb C [\mathrm{Mat}_n]_q\cdot({\rm det}_q \mathbf z)^{-k})_j=\mathbb C
[\mathrm{Mat}_n]_{q,nk+j} \cdot({\rm det}_q \mathbf z)^{-k}
\\=\bigoplus \limits_{\substack{k_n \geq 0,
\\k_1+\ldots+k_n=nk+j}} \,V_{\overline{\mathbf k}} \cdot({\rm det}_q \mathbf z)^{-k}=
\bigoplus \limits_{\substack{k_n \geq -k,
\\k_1+\ldots+k_n=j}} \,V_{\overline{\mathbf k}}.
\end{multline*}
(Since there is no zero divisors in $\mathbb C [\mathrm{Mat}_n]_{q,\det_q
\mathbf z}$, the proof of statement (\ref{l}) follows from the last
equality.)\hfill $\square$

\medskip

{\it Remark.} It can be easily verified that $v^h_{\overline{\mathbf k}}$
is a $U_q\mathfrak{k}$-highest vector and with weight
$(k_1-k_2,\ldots,k_{n-1}-k_n,2k_n+\alpha-\beta,k_{n-1}-k_n,\ldots,k_1-k_2)$.
Then the highest weight of simple $U_q\mathfrak{k}$-module
$V_{\overline{\mathbf k}}$ is equal to
$(k_1-k_2,\ldots,k_{n-1}-k_n,2k_n+\alpha-\beta,k_{n-1}-k_n,\ldots,k_1-k_2)$.

\medskip

In the classical case $ \mathfrak{sl}_{2n}=\mathfrak{p}^{-} \oplus
\mathfrak{k} \oplus \mathfrak{p}^{+},$ where
\begin{equation*}
\mathfrak{p}^{-} = \left \{\left.\begin{pmatrix}0 & 0 \\A & 0\end{pmatrix}
\right | A \in \mathrm{Mat}_{n,n}(\mathbb C)\right \}, \quad
\mathfrak{p}^{+} =\left \{\left.
\begin{pmatrix} 0 & A\\0 & 0\end{pmatrix} \right | A \in \mathrm{Mat}_{n,n}(\mathbb C)
\right \}.
\end{equation*}
Therefore $U\mathfrak{sl}_{2n} \simeq U\mathfrak{p}^{-} \otimes
U\mathfrak{k} \otimes U\mathfrak{p}^{+}$ as $U\mathfrak{k}$-modules
($U\mathfrak{p}^{-}$ and $U\mathfrak{p}^{-}$ are $U\mathfrak{k}$-modules
under the adjoint action).

In the quantum case we have an analogue of this decomposition obtained by
Jakobsen in \cite{Jak}. A quantum analogue $\ad_a, a \in U_q
\mathfrak{sl}_{2n}$ of the adjoint action is introduced via the Hopf
algebra structure of $U_q \mathfrak{sl}_{2n}$. There are $n^2$-dimensional
vector subspaces $\mathfrak{p}_q^{+} = U_q\mathfrak{k}\cdot E_n$,
$\mathfrak{p}_q^{-}=U_q \mathfrak{k}\cdot (K_nF_n)$, which are
$U_q\mathfrak{k}$-invariant under the adjoint action.\footnote{%The adjoint
%action is defined by the operators $ad_a $, $a \in U_q \mathfrak{sl}_{2n}$.
The operator $\ad_a$ is defined on $b \in U_q \mathfrak{sl}_{2n}$ in the
following way: $\ad_a(b)=\sum S(a^{'})ba^{''}$, where $\triangle a=\sum
a^{'} \otimes a^{''}$ is the comultiplication , $S$ is the antipode in $U_q
\mathfrak{sl}_{2n}$.} Instead of $U\mathfrak{p}^{-},U\mathfrak{p}^{+}$,
there are the subalgebras $U_q\mathfrak{p}^{-},U_q\mathfrak{p}^{+} \subset
U_q\mathfrak{sl}_{2n}$ generated by
$\mathfrak{p}_q^{-},\mathfrak{p}_q^{+}$, respectively. The algebras
$U_q\mathfrak{p}^{-}$ and $U_q\mathfrak{p}^{+}$ are
$U_q\mathfrak{k}$-modules under the adjoint action. Therefore in the
quantum case we get $U_q\mathfrak{sl}_{2n} \simeq U_q\mathfrak{p}^{-}
\otimes U_q\mathfrak{k} \otimes U_q\mathfrak{p}^{+}$ as
$U_q\mathfrak{k}$-modules (see \cite{Jak}). It's worthwhile to note that
$U_q\mathfrak{p}^{-}$ and $U_q\mathfrak{p}^{+}$ are not Hopf subalgebras
unlike the classical case.

%Let us recall some basic notions from the theory of Verma modules over $U_q
%\mathfrak{sl}_{2n}$ (see \cite{Jantzen},chap.5). Let $P_+=\mathbb
%Z_+^{2n-1}$ and $\lambda=(\lambda_1,\ldots,\lambda_{2n-1}) \in P_+$. The
%Verma module $M(\lambda)$ can be defined by its generator $v(\lambda)$
% and the relations
%$$
%E_iv(\lambda)=0,\qquad K_iv(\lambda)=q^{\lambda_i}v(\lambda).
%$$
%Let $K(\lambda)$ be the maximal proper submodule of $M(\lambda)$. Then
%$L(\lambda)=M(\lambda) \diagup K(\lambda)$ is a simple finite-dimensional
%weight $U_q \mathfrak{sl}_{2n}$-module with highest weight $\lambda$.

In the last part of this section we describe how each $U_q
\mathfrak{k}$-isotypic component $V_{\overline{\mathbf k}}$ transforms
under the action of $\mathfrak{p}_q^{-}$ and $\mathfrak{p}_q^{+}$. This
allows one to understand how $V$ transforms under the
$U_{q}\mathfrak{sl}_{2n}$-action. Since $U_q \mathfrak{k} =
U_q\mathfrak{s(gl}_n \times \mathfrak{gl}_n) \simeq \mathbb C[K_{0}^{\pm
1}] \otimes (U_q\mathfrak{sl}_n \otimes U_q\mathfrak{sl}_n)=\mathbb
C[K_{0}^{\pm 1}] \otimes U_q\mathfrak{k}_{ss}$, where
\begin{equation} \label {K_0}
K_0=K_1 K_2^2K_3^3 \ldots K_n^nK_{n+1}^{n-1}\ldots K_{2n-2}^2K_{2n-1},
\end{equation}
and $\pi_{\alpha,\beta}(K_0)$ acts by scalar multiplications in every
isotypic component, we see that $V_{\overline{\mathbf k}}$ is a simple
$U_q\mathfrak{sl}_n \otimes U_q\mathfrak{sl}_n$-module. Hence the
$U_q\mathfrak{sl}_n \otimes U_q\mathfrak{sl}_n$-module
$V_{\overline{\mathbf k}}$ decomposes into a tensor product of
$U_{q}\mathfrak{sl}_n$-modules: $V_{\overline{\mathbf k}} \simeq
V_{\overline{\mathbf k}}^{(1)} \otimes V_{\overline{\mathbf k}}^{(2)}$ with
$V_{\overline{\mathbf k}}^{(1)} = L(k_1-k_2,\ldots,k_{n-1}-k_n)$,
$V_{\overline{\mathbf k}}^{(2)} = L^*(k_1-k_2,\ldots,k_{n-1}-k_n)$, where
we denote by $L(k_1-k_2,\ldots,k_{n-1}-k_n)$ and
$L^*(k_1-k_2,\ldots,k_{n-1}-k_n)$ the simple finite dimensional
$U_q\mathfrak{sl}_n$-modules with highest weights
$(k_1-k_2,\ldots,k_{n-1}-k_n)$ and $(k_{n-1}-k_n,\ldots,k_1-k_2)$,
respectively.

We can equip the vector spaces $V_{\overline{\mathbf k}}^{(1)}$ and
$V_{\overline{\mathbf k}}^{(2)}$ with the structure of $U_q\mathfrak{sl}_n
\otimes U_q\mathfrak{sl}_n$-modules as follows:
\begin{equation*}
(\xi \otimes \eta)(v) = \xi \cdot (\varepsilon(\eta)v), \quad (\xi \otimes
\eta)(v^{*}) = \eta \cdot (\varepsilon(\xi)v^{*}),
\end{equation*}
for all $\xi, \eta \in U_q\mathfrak{sl}_n$, $v \in V_{\overline{\mathbf
k}}^{(1)}$, $v^{*} \in V_{\overline{\mathbf k}}^{(2)}$, where $\varepsilon$
denotes the counit of $U_{q}\mathfrak{sl}_n$. Note that as
$U_q\mathfrak{sl}_n \otimes U_q\mathfrak{sl}_n$-modules
\begin{equation*}
\mathfrak{p}_q^+ \simeq \mathbb C^n \otimes (\mathbb C^n)^{*}, \quad
\mathfrak{p}_q^- \simeq (\mathbb C^n)^{*} \otimes \mathbb C^n,
\end{equation*}
where $\mathbb C^n$ is the vector representation of $U_q\mathfrak{sl}_n$.
Consider the natural maps
\begin{equation*}
m^+_{\overline{\mathbf k}}:\mathfrak{p}_q^+ \otimes V_{\overline{\mathbf
k}} \longrightarrow V, \quad m^-_{\overline{\mathbf k}}:\mathfrak{p}_q^-
\otimes V_{\overline {\mathbf k}} \longrightarrow V.
\end{equation*}
Since there exist the $U_q\mathfrak{sl}_n \otimes
U_q\mathfrak{sl}_n$-homomorphisms
\begin{align*}
&\mathfrak{p}_q^{+} \otimes V_{\overline{\mathbf k}} \simeq \mathbb C^n
\otimes (\mathbb C^n)^{*} \otimes V_{\overline{\mathbf k}}^{(1)} \otimes
V_{\overline{\mathbf k}}^{(2)} \simeq \mathbb C^n \otimes
V_{\overline{\mathbf k}}^{(1)} \otimes(\mathbb C^n)^{*} \otimes
V_{\overline{\mathbf k}}^{(2)},
\\& \mathfrak{p}_q^{-} \otimes V_{\overline{\mathbf k}} \simeq
(\mathbb C^n)^{*} \otimes \mathbb C^n \otimes V_{\overline{\mathbf
k}}^{(1)} \otimes V_{\overline{\mathbf k}}^{(2)} \simeq (\mathbb C^n)^{*}
\otimes V_{\overline{\mathbf k}}^{(1)} \otimes \mathbb C^n \otimes
V_{\overline{\mathbf k}}^{(2)},
\end{align*}
we have the well-defined morphisms
\begin{eqnarray*}
\mathcal{M^+}_{\overline{\mathbf k}}:&\mathbb C^n \otimes
V_{\overline{\mathbf k}}^{(1)} \otimes(\mathbb C^n)^{*} \otimes
V_{\overline{\mathbf k}}^{(2)} &\longrightarrow V,
\\\mathcal{M^-}_{\overline{\mathbf k}}:&(\mathbb C^n)^{*} \otimes
V_{\overline{\mathbf k}}^{(1)} \otimes \mathbb C^n \otimes
V_{\overline{\mathbf k}}^{(2)} &\longrightarrow V.
\end{eqnarray*}

For instance, if we consider a $U_{q} \mathfrak {sl}_n$-highest vector
$v_{1} \in \mathbb C^n \otimes V_{\overline{\mathbf k}}^{(1)}$ and a $U_{q}
\mathfrak {sl}_n$-highest vector $v_{2} \in (\mathbb C^n)^{*} \otimes
V_{\overline{\mathbf k}}^{(2)}$, then $\mathcal{M^+} _{\overline{\mathbf
k}}(v_1 \otimes v_2)$ is a $U_{q} \mathfrak {sl}_n \otimes U_{q} \mathfrak
{sl}_n$-highest vector (or, equivalently, a $U_{q}\mathfrak{k}$-highest
vector) in $V$.

In the sequel we are going to get explicit formulas for $U_q \mathfrak
{sl}_n$-highest vectors $\zeta_j \in \mathbb C^n \otimes
L(k_1-k_2,\ldots,k_{n-1}-k_n)$, $j=1,\ldots,n$ with weights
$(k_1,\ldots,k_j+1,\ldots,k_n)$, respectively.

\medskip

In the classical case auxiliary elements $\mathrm{F_{mj}}$ of the universal
enveloping algebra $U \mathfrak {sl}_n$ are used in such formulas.
\begin{lemma}
(\cite{Lee}, lemma 3.4) \ \ Let $1 \leq k \leq n-1$ and $1 \leq m < j \leq
n$.
\begin{enumerate}
\item If $1 \leq k < m$ or $j < k \leq n$, then $\mathrm{E_k F_{mj}}=
\mathrm{F_{mj}E_k}.$ \item If $k=m$, then $\mathrm{E_m F_{mj}} \equiv
\mathrm{F_{m+1,j}(H_m+\ldots+H_{j-1}+j-m-1)} \;(\!\!\!\!\mod \,U \mathfrak
{sl}_n \cdot E_m).$ \item If $m < k \leq j$, then $\mathrm{E_k F_{mj}
\equiv 0 \;(\!\!\!\!\mod \,U \mathfrak {sl}_n \cdot E_k)}.$ \hfill
$\square$
\end{enumerate}
\end{lemma}
Explicit formulas for the elements $\mathrm{F_{mj}}$ are used for the proof
of lemma.
\begin{equation*}
\mathrm{F}_{mj}=\mathrm{F_{m+1,j}F_m}+\sum_{t=m+2}^j (-1)^{t+m+1}
\mathrm{F}_{tj} \,\mathrm{\ad_{F_{t-1}}\ldots \ad_{F_m} H(j;m+1,t-1)},
\end{equation*}
where $\mathrm{H(j;p,s)=\prod \limits_{a=p}^s (H_a+\ldots+H_{j-1}+j-a)}$.
\par We find quantum analogues of the previous lemma and the elements in
$U_q \mathfrak {sl}_n$.

For $1 \leq m \leq j \leq n$ define $F_{mj} \in U_q \mathfrak {sl}_n$
inductively as follows:
$$
F_{jj}=1, \quad F_{j-1,j}=F_{j-1}K_{j-1},
$$
$$
F_{mj}=F_{m+1,j}\,F_mK_m+\sum_{s=m+2}^{j} (-1)^{s+m+1}F_{sj}
\,\ad_{F_{s-1}} \ldots \ad_{F_{m+1}}(F_mK_m) K(j,m+1,s-1),
$$
where $K(j,p,r)=\prod \limits_{a=p}^{r} q^{j-a}K_a \ldots K_{j-1}
[H_a+\ldots+H_{j-1}+j-a]_q$.

Here and everywhere below we use the standard notation $[x]_q =
\frac{q^x-q^{-x}} {q-q^{-1}}$.
\begin{lemma}\label{l_1}
The following relations are satisfied:
\begin{align}
&\!\!\!1. \,K_iF_{mj}=F_{mj}K_i \quad \text{for} \; 1 \leq i< m-1 \;
\text{or} \; j<i \leq n. \label{l_1_1}
\\& \!\!\!2. \,K_jF_{mj}=qF_{mj}K_j,\quad K_{m-1}F_{mj}=qF_{mj}K_{m-1}. \label{l_1_2}
\\& \!\!\!3. \,qK_{j-1}F_{mj}=F_{mj}K_{j-1},\quad qK_{m}F_{mj}=F_{mj}K_{m}.
\label{l_1_3}
\\& \!\!\!4. \,E_iF_{mj}=F_{mj}E_i \quad \text{for} \; 1 \leq i < m-1
\; \text{or} \; j <i \leq n. \label{l_1_4}
\\& \!\!\!5. \,E_{m-1}F_{mj}= q F_{mj}E_{m-1},\quad E_jF_{mj}=qF_{mj}E_j.
\label{l_1_5}
\\& \!\!\!6. \,E_iF_{mj} \equiv 0 \;(\!\!\!\!\!\mod U_q \mathfrak {sl}_n \cdot E_i)
\quad \text{for} \; m<i<j. \label{l_1_6}
\\& \!\!\!7. \,E_mF_{mj} \!\equiv \!F_{m+1,j}\,
q^{j-m}\!K_m\!\ldots \!K_{j-1}[H_m+\ldots+H_{j-1}+j-m-1]_q
\;(\!\!\!\!\!\mod U_q \mathfrak {sl}_n \!\cdot \!E_m). \label{l_1_7}
\end{align}
\end{lemma}
This lemma is proved in Appendix.\hfill $\square$

Let $(\varepsilon_1,\ldots,\varepsilon_n)$ be the standard basis for
$\mathbb C^{n}$. Suppose $u \in L(k_1-k_2,\ldots,k_{n-1}-k_n)$ is a $U_{q}
\mathfrak {sl}_n$-highest vector with weight
$(k_1-k_2,\ldots,k_{n-1}-k_n)$.
\begin{proposition} \label{a}
Define vectors $\{\zeta_j\}_{j=1}^n$ as follows:
$$
\zeta_j=\sum \limits_{m=1}^j (-q^2)^{m-1} \varepsilon_m \otimes
F_{mj}K_-(j,1,m-1)u \;\in \;\mathbb C^n \otimes
L(k_1-k_2,\ldots,k_{n-1}-k_n),
$$
where $K_-(j,p,r)=\prod \limits_{a=p}^{r} q^{j-a-1}K_a \ldots K_{j-1}
[H_a+\ldots+H_{j-1}+j-a-1]_q$. Then $\zeta_j$ is a $U_{q} \mathfrak
{sl}_n$-highest vector\footnote{I.e. $E_i \zeta_j=0$ for all
$i=1,\ldots,n-1$.} with weight
$(k_1-k_2,\ldots,k_{j-1}-k_j-1,k_j+1-k_{j+1},\ldots,k_{n-1}-k_n)$ for
$j=1,\ldots,n$.
\end{proposition}

{\bf Proof.} Using Lemma \ref{l_1}, it is easy to prove that $\zeta_j$ are
weight vectors. We claim that $E_i \,\zeta_j=0$ for all $1 \leq i \leq
n,\,1 \leq j \leq n$. Indeed, by Lemma \ref{l_1}
\begin{multline*}
E_i \zeta_j=E_i\,(\sum \limits_{m=1}^j (-q^2)^{m-1} \varepsilon_m \otimes
F_{mj}K_-(j,1,m-1)u)
\\ =\sum \limits_{m=1}^j (-q^2)^{m-1}
E_i\,(\varepsilon_m \otimes F_{mj}\, K_-(j,1,m-1)u) =(-q^2)^i
\varepsilon_{i} \otimes F_{i+1,j}\,K_-(j,1,i)u
\\+(-q^2)^{i-1} E_i\,(\varepsilon_i \otimes F_{i,j}\,K_-(j,1,i-1)u)=
(-q^2)^{i-1}(-q^2 \varepsilon_i \otimes F_{i+1,j}\,K_-(j,1,i) u
\\+q \varepsilon_i \otimes F_{i+1,j}\, q^{j-i}K_i\ldots K_{j-1}
[H_i+\ldots+H_{j-1}+j-i-1]_q \,K_-(j,1,i-1)u)=0. \tag*{$\square$}
\end{multline*}

Similarly, we are going to get explicit formulas for $U_{q} \mathfrak
{sl}_n$-highest vectors $\xi_j \in (\mathbb C^n)^{*} \otimes
L^*(k_1-k_2,\ldots,k_{n-1}-k_n)$. For $1 \leq r \leq t \leq n$ introduce
the elements $S_{rt} \in U_{q} \mathfrak {sl}_n$ as
follows:\footnote{Classic analogues of these elements were investigated by
Lee in \cite{Lee}.}
$$
S_{tt}=1,\quad S_{t-1,t}=F_tK_t,
$$
$$
S_{rt}=S_{r,t-1}F_tK_t+\sum_{s=r+1}^{t-1}S_{r,s-1} \,\ad_{F_s}\ldots
\ad_{F_{t-1}}(F_tK_t)\, L(t,s,t-1),
$$
where $L(j,p,r)=\prod \limits_{a=p}^r q^{a-j}K_{j+1}\ldots
K_a[H_{j+1}+\ldots+H_a+a-j]_q$.
\begin{lemma} \label{l_2}
The following relations are satisfied:
\begin{enumerate}
\item $ K_iS_{rt}=S_{rt}K_i$ for $1 \leq i < r$ or $t+1 < i \leq n$. \item
$K_rS_{rt}=qK_rS_{rt}, \quad K_{t+1}S_{rt}=qK_{t+1}S_{rt}$. \item
$K_{r+1}S_{rt}=q^{-1}K_{r+1}S_{rt}, \quad K_tS_{rt}=q^{-1}K_tS_{rt}$. \item
$E_iS_{rt}=S_{rt}E_i$ for $1 \leq i < r$ or $t+1 < i \leq n$. \item
$E_rS_{rt}=qS_{rt}E_r, \quad E_{t+1}S_{r,t}=qS_{r,t}E_{t+1}$. \item
$E_iS_{rt} \equiv 0 \;(\!\!\!\!\mod \,U_{q} \mathfrak {sl}_n \cdot E_i)$
for $r<i<t$. \item $E_tS_{rt} \equiv -S_{r,t-1} q^{t-r}K_{r+1}\ldots K_t
[H_{r+1}+\ldots+H_t+t-r-1]_q \;(\!\!\!\!\mod \,U_{q} \mathfrak {sl}_n \cdot
E_t)$.
\end{enumerate}
\end{lemma}
The proof of this lemma is completely analogous to the proof of Lemma
\ref{l_1}.

Let $(\varepsilon^*_1,\ldots,\varepsilon^*_n)$ be the basis for $(\mathbb
C^n)^*$ dual to the basis $(\varepsilon_1,\ldots,\varepsilon_n)$ for
$\mathbb C^n$. Suppose that $u^* \in L^*(k_1-k_2,\ldots,k_{n-1}-k_n)$ is a
$U_{q} \mathfrak {sl}_n$-highest vector with weight
$(k_{n-1}-k_n,\ldots,k_1-k_2)$. The proof of the next statement is similar
to the proof of Proposition \ref{a}.
\begin{proposition} \label{b}
Define vectors $\{\xi_j\}_{j=1}^n$ as follows:
$$
\xi_j=\sum \limits_{m=j}^n \varepsilon_m^* \, \otimes \,
S_{jm}\,L_-(j,m+1,n)\,u^*\;\in \; (\mathbb C^n)^{*} \otimes
L^*(k_1-k_2,\ldots,k_{n-1}-k_n)
$$
where $L_-(j,p,r)=\prod \limits_{a=p}^r q^{a-j-1}K_{j+1}\ldots K_a
\,[H_{j+1}+\ldots+H_a+a-j-1]_q$. Then $\xi_j$ is a $U_{q} \mathfrak
{sl}_n$-highest vector with weight
$(k_{n-1}-k_n,\ldots,k_{n-j+1}+1-k_{n-j+2},k_{n-j}-k_{n-j+1}-1,\ldots,k_1-k_2)$
for $j=1,\ldots,n$. \hfill $\square$
\end{proposition}

It follows from Propositions \ref{a} and \ref{b} that $\mathcal{M^+}
_{\overline{\mathbf k}}:\mathbb C^n \otimes V_{\overline{\mathbf k}}^{(1)}
\otimes(\mathbb C^n)^{*} \otimes V_{\overline{\mathbf k}}^{(2)} \rightarrow
\bigoplus_{j=1}^n V_{\bf \overline{k} +e_j}.$ For all $j,k=1,\ldots,n$ the
vectors $\mathcal{M^+}_{\overline{\mathbf k}} (\zeta_j \otimes \xi_k)$ are
$U_q \mathfrak{sl}_n \otimes U_q \mathfrak{sl}_n$-highest vectors in $V$.
By the action of $\pi_{\alpha,\beta}(K_0)$ (see (\ref{K_0})), the vector
$\mathcal{M^+} _{\overline{\mathbf k}} (\zeta_j \otimes \xi_k)$ is a $U_q
\mathfrak{k}$-highest vector in $V_{\bf \overline{k}+e_j}$ if and only if
$k=n-j+1$. Since every isotypic components occurs with multiplicity one,
$\mathcal{M^+}_{\overline{\mathbf k}} (\zeta_j \otimes \xi_{n-j+1})=c_j
\cdot v^h_{\bf \overline{k}+e_1} = c_j \cdot (z^{\wedge
1})^{k_1-k_2+1}\ldots(z^{\wedge n-1})^{k_{n-1}-k_n} (z^{\wedge n})^{k_n}$
for some $c_j \in \mathbb C$. (Here and below we suppose that if
$\overline{\mathbf m}=(m_1,\ldots,m_n) \not \in \widehat{K}$, then
$V_{\overline{\mathbf m}}=0$ and $v^h_{\overline{\mathbf m}}=0$.)

The proof of the next statement, reduced to computation of $c_j$, is given
in Appendix.

\begin{proposition}\label{pr_2}
For every $j=1,\ldots,n,\,\overline{\mathbf k} \in \widehat{K}$
$$
\mathcal{M^+}_{\overline{\mathbf k}}(\zeta_j \otimes \xi_{n-j+1}) =
c_j(\beta,k_j) v^h_{\bf \overline{k}+e_j},
$$
where $c_j(\beta,k_j) = q^{-\beta-n/2}[\beta-k_j+j-1]_q \omega_j(\overline
{\mathbf k},q)$ and
 $\omega_j(\overline {\mathbf k},q) \not =0$ for all $\overline{\mathbf k}
\in \widehat{K}$.
\end{proposition}

We deduce sufficient conditions for reducibility of $\pi_{\alpha,\beta}$
from Proposition \ref{pr_2}.

Let $\alpha,\beta$ be fixed. For any $j=1,\ldots,n$ and $\overline{\mathbf
k} \in \widehat{K}$ if $c_j(\beta,k_j) \neq 0$, then there exist $v \in
V_{\overline{\mathbf k}}, \; \xi \in \mathfrak{p}_q^{+}$ such that
$\pi_{\alpha,\,\beta}(\xi) \cdot v \in V_{\bf \overline{k} +e_j}$. That
means $\pi_{\alpha,\,\beta}(U_{q}\mathfrak{sl}_{2n}) \cdot
V_{\overline{\mathbf k}} \supset V_{\bf \overline{k} +e_j}$.

Let us consider in details other cases, i.e., let $c_j(\beta,k_j)=0$ for
some $k_j$. For fixed $\beta$, by Proposition \ref{pr_2} and (\ref{lim}),
the equation $c_j(\beta,k_j) = 0$ is equivalent to $\beta-k_j+j-1=0$.
\begin{corollary}\label{cor_2}
For all $j=1,\ldots,n,\,\overline{\mathbf k} \in \widehat{K}$, the subspace
$V^j_{\leq k} \stackrel {\rm def}{=} \bigoplus \limits_{\{\overline{\mathbf
k'} \in \widehat{K}| k \geq k'_j\}} V_{\overline{\mathbf k'}}$ is a
$U_q\mathfrak{sl}_{2n}$-submodule in $V$ iff $\beta-k+j-1=0$.
\end{corollary}

{\bf Proof.} Let $j=1$, the other cases are similar. The necessity easily
follows from the above. Prove the sufficiency. If $\beta-k_1=0$, then
$\mathcal{M^+}_{\overline{\mathbf k}}(\mathfrak{p}_q^{+} \otimes
V_{\overline{\mathbf k}}) \subset \bigoplus \limits_{j=2}^n V_{\bf
\overline {k} +e_j}$. Introduce the natural filtration on
$U_q\mathfrak{p}^{+}$ (here $U_q\mathfrak{p}^{+}$ is the algebra generated
by $\mathfrak{p}_q^{+}$) in the following way: $U_q\mathfrak{p}^{+}=\bigcup
\limits_{n \geq 0} (U_q\mathfrak{p}^{+})^{(n)}$. Then
$\pi_{\alpha,\,\beta}((U_q\mathfrak{p}^{+})^{(1)})(V_{\overline{\mathbf
k}})\subset V_{\overline{\mathbf k}} \oplus (\bigoplus \limits_{j=2}^n
V_{\bf \overline{k} +e_j})$. In the same way,
$\pi_{\alpha,\,\beta}((U_q\mathfrak{p}^{+})^{(2)})(V_{\overline{\mathbf
k}}) \subset
\pi_{\alpha,\,\beta}((U_q\mathfrak{p}^{+})^{(1)})(V_{\overline{\mathbf k}}
\oplus (\bigoplus \limits_{j=2}^n V_{\bf \overline{k}+e_j})) \subset
V_{\overline{\mathbf k}} \oplus (\bigoplus \limits_{j=2}^n V_{\bf
\overline{k}+e_j}) \oplus (\bigoplus \limits_{n \geq j_1\geq j_2\geq 2}
V_{\bf \overline{k}+e_{j_1}+e_{j_2}})$. Then
$\pi_{\alpha,\,\beta}(U_q\mathfrak{p}^{+})(V_{\overline{\mathbf k}})
\subset \bigoplus \limits_{m=0}^{\infty}(\bigoplus \limits_{n \geq j_1\geq
\ldots \geq j_m\geq 2}^n V_{\bf \overline{k}+e_{j_1}+\ldots+e_{j_m}})$, and
\begin{multline*}
\hspace{-10pt}
\pi_{\alpha,\,\beta}(U_q\mathfrak{sl}_{2n})(V_{\overline{\mathbf k}})
\subset \pi_{\alpha,\,\beta}(U_q\mathfrak{p}^{-})
\pi_{\alpha,\,\beta}(U_q\mathfrak{k})
(\pi_{\alpha,\,\beta}(U_q\mathfrak{p}^{+}) V_{\overline{\mathbf k}})
\\ \subset \pi_{\alpha,\,\beta}(U_q\mathfrak{p}^{-}) \pi_{\alpha,\,\beta}(U_q
\mathfrak{k}) (\bigoplus_{m=0}^{\infty} (\bigoplus \limits_{n \geq j_1\geq
\ldots \geq j_m\geq 2} V_{\bf \overline {k}+e_{j_1}+\ldots+e_{j_m}}))
\\ \subset \pi_{\alpha,\,\beta}(U_q\mathfrak{p}^{-})(\bigoplus_{m \geq 0} (\bigoplus
\limits_{n \geq j_1\geq \ldots \geq j_m\geq 2} V_{\bf \overline{k}
+e_{j_1}+\ldots+e_{j_m}})) \subset V^1_{\leq k}.
\end{multline*}
Obviously the subspace $V^1_{\leq k}$ is a
$U_q\mathfrak{sl}_{2n}$-submodule in $V$. \hfill $\square$

\medskip

By the same arguments as in Propositions \ref{a}, \ref{b} and \ref{pr_2},
one has
\begin{proposition} \label{c}
Define vectors $\{\xi'_j\}_{j=1}^n$ as follows:
$$
\xi'_j=\sum \limits_{m=j}^n \varepsilon_m^* \otimes S_{jm}L_-(j,m+1,n)u
\;\in \; (\mathbb C^n)^{*} \otimes L(k_1-k_2,\ldots,k_{n-1}-k_n)
$$
where $L_-(j,p,r)=\prod \limits_{a=p}^r q^{a-j-1}K_{j+1}\ldots K_a
[H_{j+1}+\ldots+H_a+a-j-1]_q$. Then $\xi'_j$ is a $U_{q} \mathfrak
{sl}_n$-highest vector with weight
$(k_1-k_2,\ldots,k_{j-1}-k_j+1,k_j-1-k_{j+1},\ldots,k_{n-1}-k_n)$ for
$j=1,\ldots,n$. \hfill $\square$
\end{proposition}
\begin{proposition}\label{d}
Define vectors $\{\zeta'_j\}_{j=1}^n$ as follows
$$
\zeta'_j=\sum \limits_{m=1}^j (-q^2)^{m-1} \varepsilon_m \otimes
F_{mj}K_-(j,1,m-1)u^* \;\in \; \mathbb C^n \otimes
L^*(k_1-k_2,\ldots,k_{n-1}-k_n)
$$
where $K_-(j,p,r)=\prod \limits_{a=p}^{r} q^{j-a-1}K_a\ldots K_{j-1}
[H_a+\ldots+H_{j-1}+j-a-1]_q$. Then $\zeta'_j$ is a $U_{q} \mathfrak
{sl}_n$-highest vector with weight
$(k_{n-1}-k_n,\ldots,k_{n-j+1}-k_{n-j}-1,k_{n-j}+1-k_{n-j-1},\ldots,k_1-k_2)$
for $j=1,\ldots,n$. \hfill $\square$
\end{proposition}

The proof of the next statement, reduced as for Proposition \ref{pr_2} to
computation of $d_j$, is given in Appendix.
\begin{proposition}\label{pr_3}
For every $j=1,\ldots,n,\,\overline{\mathbf k} \in \widehat{K}$
$$
\mathcal{M^-}_{\overline{\mathbf k}}(\xi'_j \otimes \zeta'_{n-j+1}) =
d_j(\alpha,k_j) v^h_{\bf \overline{k}-e_j},
$$
where $d_j(\alpha,k_j) = q^{\alpha+n/2}[\alpha+k_j+n-j]_q
\varpi_j(\overline {\mathbf k},q)$ and $\varpi_j(\overline {\mathbf k},q)
\not=0$ for all $\overline{\mathbf k} \in \widehat{K}$. \hfill $\square$
\end{proposition}

By (\ref{lim}) and Proposition \ref{pr_3}, we see that the equations
$d_j(\alpha,k_j) =0$ and $\alpha+k_j+n-j=0$ are equivalent.
\begin{corollary}\label{cor_3}
For all $j=1,\ldots,n,\,\overline{\mathbf k} \in \widehat{K}$ the subspace
$V^j_{\geq k} \stackrel {\rm def}{=} \bigoplus \limits_{\{\overline{\mathbf
k'} \in \widehat{K}| k'_j \geq k\}} V_{\overline{\mathbf k'}}$ is a
$U_q\mathfrak{sl}_{2n}$-submodule in $V$ iff $\alpha+k_j+n-j=0$. \hfill
$\square$
\end{corollary}

\section{Reducibility of $\pi_{\alpha,\beta}$}\label{reduce}
\begin{proposition} \label{pr_rez}
The representation $\pi_{\alpha,\beta}$ is irreducible if and only if
$\alpha,\beta$ satisfy the following equivalent conditions:\footnote{Since
$\alpha-\beta \in \mathbb Z$, these conditions are equivalent.}
\begin{equation*}
1.\,\,\alpha \not \in \mathbb Z; \qquad 2.\,\,\beta \not \in \mathbb Z.
\end{equation*}
\end{proposition}

{\bf Proof.} Suppose $\alpha \not \in \mathbb Z$, $\beta \not \in \mathbb
Z$. Consider the system of equations
\begin{equation*}
\begin{cases}
\beta-k_1=0,
\\\beta-k_2+1=0,
\\\ldots\ldots\ldots\ldots..
\\\beta-k_n+n-1=0,
\\\alpha+k_1+n-1=0,
\\\ldots\ldots\ldots\ldots\ldots
\\\alpha+k_n=0.
\end{cases}
\end{equation*}
This system has no integral solution. Therefore $c_j(\beta,k_j)$ and
$d_j(\alpha,k_j)$ do not vanish. Let $W$ be a
$U_q\mathfrak{sl}_{2n}$-submodule of $V$. Then $W=\bigoplus
\limits_{\overline {\mathbf k} \in I} V_{\overline {\mathbf k}}$ for some
$I \subset \widehat{K}$. Then, for all $\overline {\mathbf k} \in I$ and
$j=1,\ldots,n$, it follows that ${\bf \overline {k} + e_j,\overline {k} -
e_j} \in I$ (if the respective indexes belong to $\widehat{K}$). Therefore
if $I \neq \emptyset$, then $I=\widehat{K}$, and the module $V$ have no
proper submodules, i.e. it is simple. Conversely, by Corollaries
\ref{cor_2} and \ref{cor_3}, if $\pi_{\alpha,\beta}$ is irreducible, then
$\alpha \not \in \mathbb Z$, $\beta \not \in \mathbb Z$. \hfill $\square$

\begin{corollary}
Let $\alpha,\beta \in \mathbb Z$, and let $W$ be the representation space
of a subrepresentation of $\pi_{\alpha,\,\beta}$. Then $W$ is a finite
intersection of some of the $U_q\mathfrak{sl}_{2n}$-modules $V^j_{\geq k},
\,V^j_{\leq k}$ defined in Corollaries \ref{cor_2},\ref{cor_3}.
\end{corollary}
The proof follows directly from the previous proof. \hfill $\square$

\medskip

Now suppose that $\alpha,\,\beta \in \mathbb Z$. We will investigate
reducibility and proper subrepresentations of $\pi_{\alpha,\,\beta}$. We
use figures as in \cite{Howe,Lee} for description.

Each $U_q\mathfrak{k}$-isotypic component $V_{\overline{\mathbf k}}$ is
assigned to the point $(k_1,\ldots,k_n) \in \mathbb R^n$. Thus
$\widehat{K}$ is assigned to the set $\mathbf
K^{+}=\{(k_1,\ldots,k_n)\,|\,k_1 \geq \ldots \geq k_n \} \subset \mathbb
R^n$. Consider $2n$ hyperplanes:
$$\mathcal{L}_j^{+}:k_j=\beta+j-1; \qquad \mathcal{L}_j^{-}: k_j=-\alpha-n+j.$$

These hyperplanes are parallel to the coordinate axis and pass through
points with integral coordinates. The distance between $\mathcal{L}_j^{+}$
and $\mathcal{L}_j^{-}$ is equal to $\alpha+\beta+n-1$.

By Corollaries \ref{cor_2} and \ref{cor_3},
\begin{align*}
\overline{\mathbf k} \in \mathcal{L}_j^{+} \quad \text{iff} \quad
U_q\mathfrak{sl}_{2n} \cdot V_{\overline{\mathbf k}} \not \supset V_{\bf
\overline{k}+e_j}; \quad \overline{\mathbf k} \in \mathcal{L}_j^{-} \quad
\text{iff} \quad U_q\mathfrak{sl}_{2n} \cdot V_{\overline{\mathbf k}} \not
\supset V_{\bf \overline{k}-e_j}.
\end{align*}

Investigate the example $n=2$. In this case $\mathcal{L}_j^{\pm},\,j=1,2$
are just lines on the plane $\mathbb{R}^2$, parallel to the coordinate
axis. Let us consider different values of $\alpha+\beta$.

Case 1. $\alpha+\beta \geq 0$. In this case the line $\mathcal{L}_1^{+}$
lies to the right of $\mathcal{L}_1^{-}$, $\mathcal{L}_2^{+}$ lies higher
than $\mathcal{L}_2^{-}$. The lines
$\mathcal{L}_1^{\pm},\mathcal{L}_2^{\pm}$ are shown on Fig.1. The
intersection point of $\mathcal{L}_1^{+}$ and $\mathcal{L}_2^{-}$ has the
coordinates $(\beta,\,-\alpha)$ and belongs to $\mathbf K^+$. Arrows
attached to $\mathcal{L}_j^{\pm}$ show the direction of isotypic components
"movement" under $\pi_{\alpha,\beta}$. There exists a unique simple
submodule $V^s=\bigoplus \limits_{\{\overline{\mathbf k} \in
\widehat{K}|k_1 \leq \beta,\,k_2 \geq -\alpha\}} V_{\overline{\mathbf k}}$
in $V$.

Case 2. $\alpha+\beta=-1$.
 In this case the lines $\mathcal{L}_1^{+}$ and $\mathcal{L}_1^{-}$,
 $\mathcal{L}_2^{+}$ and $\mathcal{L}_2^{-}$ coincide. The intersection
point of the lines $\mathcal{L}_1^{+}$ and $\mathcal{L}_2^{+}$ does not
belong to $\mathbf K^+$ (Fig.2). There are two simple submodules in $V$:
$V^s_1= \bigoplus \limits_{\{\overline{\mathbf k} \in \widehat{K}| k_1=
-1-\alpha\}}V_{\overline{\mathbf k}}$ and $V^s_2= \bigoplus \limits_{
\{\overline{\mathbf k} \in \widehat{K}| k_2 =-\alpha\}}
V_{\overline{\mathbf k}}$.

Case 3. $\alpha+\beta=-2$. In this case the line $\mathcal{L}_1^{+}$ lies
to the left of $\mathcal{L}_1^{-}$, $\mathcal{L}_2^{+}$ lies lower than
$\mathcal{L}_2^{-}$. However, the lines $\mathcal{L}_1^{-}$ and
$\mathcal{L}_2^{+}$ intersect in the point with coordinates
$(-\alpha-1,\,\beta+1)$ (see Fig.3). Besides, the distance between
$\mathcal{L}_j^{+}$ and $\mathcal{L}_j^{-}$is equal to 1. This shows that
$V$ is a direct sum of three submodules:
$$
V^s_1= \bigoplus \limits_{\{\overline{\mathbf k} \in \widehat{K}|k_1 \leq
\beta\}} V_{\overline{\mathbf k}}, \quad V^s_2= \bigoplus
\limits_{\{\overline{\mathbf k} \in \widehat{K}|k_2 \geq -\alpha\}}
V_{\overline{\mathbf k}},\quad V^s_3= \bigoplus
\limits_{\{\overline{\mathbf k} \in \widehat{K}|k_1 \geq -\alpha-1,k_2 \leq
\beta+1\}} V_{\overline{\mathbf k}}.
$$

Case 4. $\alpha+\beta \leq -3$. In this case the intersection point of
$\mathcal{L}_1^{+}$ and $\mathcal{L}_1^{+}$ belongs to $\mathbf K^+$ (see
Fig.4). Also, there are simple submodules $V^s_1$, $V^s_2$, $V^s_3$ in $V$,
but $V$ does not decompose into their direct sum.

\medskip

Turn now to the general case. Consider all possible values of $\alpha+\beta
+n-1$.

{\bf Case 1.} $\alpha+\beta +n-1\geq 1$. In this case the hyperplanes
$\mathcal{L}_j^{\pm},\,j=1,\ldots,n$ bound in $\mathbf K^+$ the subset that
corresponds to a unique simple {\it finite dimensional} submodule
$$V^s=\bigoplus \limits_{\{\overline{\mathbf k} \in \widehat{K}|-\alpha-n+j
\leq k_j \leq \beta+j-1 \text{ for all } j=1, \ldots , n\}}
V_{\overline{\mathbf k}}.$$

{\bf Case 2.} $\alpha+\beta+n-1=0$. In this case the hyperplanes
$\mathcal{L}_j^{+}$ and $\mathcal{L}_j^{-}$ coincide. There are $n$ simple
submodules in $V$:
\begin{equation}\label{c_2_submod}
V^s_j=\bigoplus \limits_{\{\overline{\mathbf k} \in
\widehat{K}|k_j=\beta+j-1\}} V_{\overline{\mathbf k}}, \quad j=1,\ldots,n.
\end{equation}

{\bf Case 3.} $\alpha+\beta=-n$. Here the distance between
$\mathcal{L}_j^{+}$ and $\mathcal{L}_j^{-}$ is equal to 1. This allows one
to decompose the set $\widehat{K}$ into a direct sum of $n+1$ subsets
$\widehat{K}_i,i=1,\ldots,n+1$, those correspond to the simple submodules:
$V^s_i= \bigoplus \limits_{\{\overline{\mathbf k} \in \widehat{K}_i\}}
V_{\overline{\mathbf k}} \;\subset V$. The subsets $\widehat{K}_i$ are
defined as follows:
$$
\widehat{K}_i=\{\overline{\mathbf k} \in \widehat{K}|k_{i-1} \geq
-\alpha-n+i-1, \beta+i-1 \geq k_i\}
$$
(for $i=1$ and $i=n+1$ we put respectively
$\widehat{K}_1=\{\overline{\mathbf k} \in \widehat{K}|k_1 \leq \beta\}$
and $\widehat{K}_{n+1}=\{\overline{\mathbf k} \in \widehat{K}|k_n \geq
-\alpha\}$).

\medskip

{\it Remark.} Since $k_j \leq \beta+j-1$ and $k_j \geq k_l$ for all $j
\leq l \leq n$, we see that $k_l \leq \beta+l-1$. By the same reason,
since $k_j \leq -\alpha-n+j$ and $k_j \geq k_l$ for all $j \geq l \geq 1$,
we see that $k_l \leq -\alpha-n+l$.

\medskip

{\bf Case 4.} $\alpha+\beta+n-1 \leq -2$. Also, there are simple submodules
corresponded to the subsets $\widehat{K}_i$. However, $V$ is not equal to
their direct sum.

Thus we have proved the following
\begin{corollary}\label{cor_4}
For $\alpha,\beta \in \mathbb Z$ the only one from the representations
$\pi_{\alpha,\beta}$ and $\pi_{-n-\beta,-n-\alpha}$ has an irreducible
finite dimensional subrepresentation.\hfill $\square$
\end{corollary}

\section {Intertwining operators} \label{eqv}

In this section we construct the intertwining operators between the
representations $\pi_{\alpha,\beta}$ and $\pi_{-n-\beta,-n-\alpha}$ for
non-integral $\alpha,\beta$. This allows one to prove Proposition
\ref{pr_eqv}.

Let $A:V \rightarrow V$ be an intertwining operator, i.e., for all $\xi \in
U_q \mathfrak{sl}_{2n},\, v \in V$, we have $A \pi_{\alpha,\beta}(\xi)(v)
=\pi_{-n-\beta,-n-\alpha}(\xi)(Av).$ The operators $\pi_{\alpha,\beta}(U_q
\mathfrak{k}_{ss})$ are independent of $\alpha,\beta$ and
$\pi_{\alpha,\beta}(K_n)=\pi_{-n-\beta,-n-\alpha}(K_n)$. Also,
$V_{\overline{\mathbf k}}$ and $V_{\overline{\mathbf m}}$ are
non-isomorphic $U_q \mathfrak{k}$-modules for $\overline{\mathbf k} \ne
\overline{\mathbf m}$. Then $A(\alpha,\beta)|_{V_{\overline{\mathbf k}}} =
a_{\overline{\mathbf k}} (\alpha,\beta)$, $a_{\overline{\mathbf k}}
(\alpha,\beta) \in \mathbb{C}$. Let us find necessary conditions for $A$ to
be an intertwining operator in terms of $a_{\overline{\mathbf k}}
(\alpha,\beta)$. By Propositions \ref{a}, \ref{b}, \ref{c}, and \ref{d}, it
follows that for all $\overline{\mathbf k} \in \widehat{K}$ there exist
$\vartheta_j,\eta_j \in U_q \mathfrak{sl}_{2n},\,j=1,\ldots,n$, such that
$\pi_{\alpha,\beta}(\eta_j)(v^h_{\overline{\mathbf k}}) = c_j(\beta,k_j)
v^h_{\bf \overline{k}+e_j}$ and
$\pi_{\alpha,\beta}(\vartheta_j)(v^h_{\overline{\mathbf k}})=
d_j(\alpha,k_j) v^h_{\bf \overline{k}-e_j}$. (Recall that
$v^h_{\overline{\mathbf k}}$ is the $U_q \mathfrak{k}$-highest vector in
$V_{\overline{\mathbf k}}$.) Therefore the necessary conditions look as
follows: for all $j=1,\ldots,n, \overline{\mathbf k} \in \widehat{K}$,
\[A \pi_{\alpha,\beta}(\eta_j)(v^h_{\overline{\mathbf k}})=
\pi_{-n-\beta,-n-\alpha}(\eta_j)(Av^h_{\overline{\mathbf k}})
\;\;\text{and} \;\;A \pi_{\alpha,\beta}(\vartheta_j)(v^h_{\overline{\mathbf
k}})=\pi_{-n-\beta,-n-\alpha}(\vartheta_j)(Av^h_{\overline{\mathbf k}}).
\]
\par Equivalently, in terms of $a_{\overline{\mathbf k}}$,
\[
a_{\bf \overline{k}+e_j} (\alpha,\beta) c_j(\beta,k_j) v^h_{\bf
\overline{k}+e_j}= a_{\overline{\mathbf k}}
(\alpha,\beta)c_j(-n-\alpha,k_j)v^h_{\bf \overline{k}+e_j},
\]
\[
a_{\bf \overline{k}-e_j} (\alpha,\beta) d_j(\alpha,k_j) v^h_{\bf
\overline{k}-e_j}= a_{\overline{\mathbf k}}
(\alpha,\beta)d_j(-n-\beta,k_j)v^h_{\bf \overline{k}-e_j}.
\]
Thus the coefficients $a_{\overline{\mathbf k}}$ of the intertwining
operator $A$ must satisfy the following conditions: for all
$j=1,\ldots,n,\overline{\mathbf k} \in \widehat{K}$,
\begin{equation*}
\frac{a_{\bf \overline{k}+e_j} (\alpha,\beta)} {a_{\overline{\mathbf
k}}(\alpha,\beta)}=\frac{c_j(-n-\alpha,k_j)}{c_j(\beta,k_j)}, \quad \quad
\frac{a_{\bf \overline{k}-e_j} (\alpha,\beta)} {a_{\overline{\mathbf
k}}(\alpha,\beta)} =\frac{d_j(-n-\beta,k_j)}{d_j(\alpha,k_j)}.
\end{equation*}
We get from Propositions \ref{pr_2} and \ref{pr_3} that for all
$j=1,\ldots,n, \overline{\mathbf k} \in \widehat{K}$,
\begin{equation*}
\frac{a_{\bf \overline{k}+e_j} (\alpha,\beta)}{a_{\overline{\mathbf k}}
(\alpha,\beta)}=q^{n+\alpha+\beta}\frac{[-n-\alpha-k_j+j-1]_q}{[\beta-k_j+j-1]_q},
\quad \frac{a_{\bf \overline{k}-e_j} (\alpha,\beta)}{a_{\overline{\mathbf
k}}(\alpha,\beta)}
=q^{-n-\beta-\alpha}\frac{[-\beta+k_j-j]_q}{[\alpha+k_j+n-j]_q}.
\end{equation*}

As we see, the coefficients $a_{\overline{\mathbf k}} (\alpha,\beta)$ are
defined up to a scalar multiplier. By additional assumption $a_{\overline
{\mathbf 0}} (\alpha,\beta)=1$, we get the explicit formulas for the
coefficients $a_{\overline{\mathbf
k}}(\alpha,\beta)=A(\alpha,\beta)|V_{\overline{\mathbf k}}$ of the
intertwining operator $A$
\begin{equation}\label{inter_op}
a_{\overline{\mathbf
k}}(\alpha,\beta) = \prod_{j=1}^n P_j(\alpha,\beta),
\end{equation}
where \begin{equation*} P_j(\alpha,\beta) =
\begin{cases}
\prod
\limits_{i=0}^{k_j-1}\frac{1-q^{2(\alpha+n+i-j+1)}}{1-q^{2(-\beta+i-j+1)}},
&\text{for} \quad k_j >0,
\\ 1, & \text{for} \quad k_j=0,
\\ \prod \limits_{i=1-k_j}^{0} \frac{1-q^{2(-\beta+i-j)}}{1-q^{2(\alpha+n+i-j)}}, &
\text{for} \quad k_j < 0.
\end{cases}
\end{equation*}

For fixed $\alpha-\beta \in \mathbb Z$, the operator $A$ is a meromorphic
operator-function with simple poles in integral points.

\section{Unitarizable representations of the degenerate principal
series}\label{unit}

In this section we find necessary and sufficient conditions for modules of
degenerate principal series and their simple submodules to be unitarizable.

Equip $U_q \mathfrak{sl}_{2n}$ with the involution $*$ as follows:
$$
E_n^*=-K_nF_n,\quad F_n^*=-E_nK_n^{-1},\quad K_n^*=K_n,
$$
$$
E_j^*=K_jF_j,\quad F_j^*=E_jK_j^{-1},\quad K_j^*=K_j,\quad j=1,\ldots,2n-1,
j \not=n.
$$
The $*$-Hopf algebra $U_q\mathfrak{su}_{n,n} \stackrel {\rm def}{=} (U_q
\mathfrak{sl}_{2n},*)$ is a q-analogue of $U \mathfrak{su}_{n,n}$, and its
subalgebra $U_q \mathfrak{s}(\mathfrak{u}_n \times \mathfrak{u}_n)
\stackrel {\rm def}{=}(U_q \mathfrak{k},*)$ is a q-analogue of $U
\mathfrak{s}(\mathfrak{u}_n \times \mathfrak{u}_n)$.

Let us introduce two auxiliary $*$-algebras $\mathrm{Pol}(S(\mathbb{U}))_q$
and $\mathrm{Pol}(\widehat{S(\mathbb{U})})_q$ (a quantum analogue of the
Shilov boundary $S(\mathbb{U})$ of the matrix ball is introduced in
\cite{Cauhy_Sege}). Equip the algebra $\mathbb{C}[\mathrm{Mat}_n]_{q,\det_q
\mathbf z}$ with the involution $*$ defined by the formula
\begin{equation*}
(z_a^b)^*=(-q)^{a+b-2n}({\rm det}_q \mathbf z)^{-1} {\rm det}_q \mathbf
z_a^b ,
\end{equation*}
where $\det_q \mathbf z_a^b$ is the q-determinant of the matrix derived
from $\mathbf z$ by deleting the line $b$ and the column $a$. Put
$\mathrm{Pol}(S(\mathbb{U}))_q = (\mathbb{C}[\mathrm{Mat}_n]_{q,\det_q
\mathbf z},*)$ and equip it with the natural structure of a $*$-module
algebra over $U_q \mathfrak{su}_{n,n}$. The involutions in
$\mathrm{Pol}(S(\mathbb{U}))_q$ and $U_q \mathfrak{su}_{n,n}$ are
compatible, i.e., for all $f \in \mathrm{Pol}(S(\mathbb{U}))_q,\,\xi \in
U_q \mathfrak{su}_{n,n}$ we have
$$
(\xi f)^*=(S( \xi ))^*f^*,
$$
where $S$ is the antipode in the Hopf algebra $U_q \mathfrak{sl}_{2n}$.

The $*$-algebra $\mathrm{Pol}(\widehat{S(\mathbb{U})})_q$ is generated by
$z_a^b$, $a,b=1,\ldots,n$, $(\det_q \mathbf z)^{-1}$, $t$ and $t^{-1}$. The
relations between $z_a^b$ and $(\det_q \mathbf z)^{-1}$ are inherited from
the $*$-algebra $\mathrm{Pol}(S(\mathbb{U}))_q$, the other relations are
provided by the following:
$$
t^{-1}t=tt^{-1}=1,\quad tt^*=t^*t,\quad tz_a^b=q^{-1}z_a^bt,\quad
t^*z_a^b=qz_a^bt^*,\quad a,b=1,\ldots,n.
$$
Consider an embedding of algebras $\mathrm{Pol}(\widehat{S(\mathbb{U})})_q
\hookrightarrow \mathbb{C}[\mathrm{Pl}_{n,2n}]_{q,t}$ which maps $t$ to $t$
and $z_a^b$ to $t^{-1} t^{\wedge n}_{\{1,\ldots,n\}J_{a\,b}}$ (see
(\ref{emb})). Using this embedding, we can extend the $U_q
\mathfrak{su}_{n,n}$-module structure from $\mathrm{Pol} (S(\mathbb{U}))_q$
onto $\mathrm{Pol}(\widehat{S(\mathbb{U})})_q$,

In \cite{Cauhy_Sege}, the invariant integral over the Shilov boundary of
the quantum matrix ball $f \mapsto \int \limits_{S(\mathbb{U})_q} f d\mu$
is defined and the following statement is actually proved.
\begin{proposition}
The linear subspace $(t^{-n})^*\cdot \mathrm{Pol}(S(\mathbb{U}))_q \cdot
t^{-n} \subset \mathrm{Pol}(\widehat{S(\mathbb{U})})_q$ is a $U_q
\mathfrak{su}_{n,n}$-module. The linear functional
$$
(t^{-n})^*\cdot f \cdot t^{-n} \mapsto \int \limits_{S(\mathbb{U})_q}fd
\mu
$$ is a $U_q \mathfrak{su}_{n,n}$-invariant integral.
\end{proposition}

The precise meaning of two next propositions will be given if we continue
$\mathrm{Pol}(\widehat{S(\mathbb{U})})_q$ via adding to the list of
generators $t^{\lambda}$, $(t^*)^{\lambda}$, $(\det_q \mathbf z)^{\lambda}$
for all $\lambda\in\mathbb{C}$. The relations between the ''new''
generators and the action of $E_j$, $F_j$, $K_j^{\pm 1}$, $j=1,\ldots,2n-1$
can be derived from the corresponding formulas for $t^m$, $(\det_q \mathbf
z)^m$ and $(t^*)^m$, where $m \in \mathbb Z$. From the previous proposition
it follows

\begin{proposition} (cf. \cite{Geom_real}, lemma 3.2) \ \
Let $Re \lambda= -n$. Then the linear subspace
$$
(({\rm det}_q \mathbf z)^{\lambda/2}t^{\lambda})^*\cdot
\mathrm{Pol}(S(\mathbb{U}))_q \cdot ({\rm det}_q \mathbf z)^{\lambda/2}
t^{\lambda} \in \mathrm{Pol}(\widehat{S(\mathbb{U})})_q
$$
is a $U_q \mathfrak{su}_{n,n}$-module. The linear functional
$$
(({\rm det}_q \mathbf z)^{\lambda/2} t^{\lambda})^*\cdot f \cdot ({\rm
det}_q \mathbf z)^{\lambda/2} t^{\lambda} \mapsto \int
\limits_{S(\mathbb{U})_q}fd \mu
$$
is a $U_q \mathfrak{su}_{n,n}$-invariant integral.
\end{proposition}

For each $\alpha,\beta \in \mathbb Z$ define an embedding
$i_{\alpha,\beta}:\, V=\mathbb{C}[\mathrm{Mat}_n]_{q,\det_q \mathbf z}
\hookrightarrow \mathrm{Pol}(\widehat{S(\mathbb{U})})_q$ by the formula
$i_{\alpha,\beta}(f)= f\cdot(\det_q \mathbf z)^{\alpha} \cdot
t^{\alpha+\beta}$ for all $f \in \mathbb{C}[\mathrm{Mat}_n]_{q,\det_q
\mathbf z}$. Using these embeddings and the commutative relations between
$t,\,t^{-1}$ and $\det_q \mathbf z$, we get
\begin{corollary}
 Let $Re(\alpha+\beta)=-n$. Then the sesquilinear form $V \times V
\rightarrow \mathbb{C}$ defined by
\begin{equation*}
<f_1,f_2> = \int \limits_{S(\mathbb{U})_q}f_2^{*}f_1 d \mu
\end{equation*}
satisfies the condition $(\pi_{\alpha,\beta}(\xi)
u,v)=(u,\pi_{\alpha,\beta}(\xi^{*})v)$ for all $u,v \in V,\, \xi \in
U_q\mathfrak{sl}_{2n}$.
\end{corollary}

Recall the definition of unitarizable module. Let A be a $*$-Hopf algebra,
$W$ an $A$-module. Then an $A$-module $W$ is {\it unitarizable} if there
exists an Hermitian form\footnote{I.e., sesquilinear Hermitian-symmetric
positive definite form.} $(\cdot,\cdot)$, which is $A$-invariant, i.e.,
\begin{equation*}
(a u,v)=(u,a^{*}v) \quad \text{for any} \quad u,v \in W,\, a \in A.
\end{equation*}

Therefore the representation $\pi_{\alpha,\beta}$ is unitary if
$Re(\alpha+\beta)=-n$. Such representations form {\it the principal unitary
series.}

\medskip

Now we are going to find all unitarizable simple modules of degenerate
principal series and their unitarizable submodules.

Weight subspaces are pairwise orthogonal with respect to every $U_q
\mathfrak{su}_{n,n}$-invariant scalar product. Therefore the isotypic
components $V_{\overline{\mathbf k}}$ are pairwise orthogonal too. From
Proposition \ref{K_reduce} and the Burnside theorem (see \cite{Kert}, \S
27), it follows that in every component $V_{\overline{\mathbf k}}$ there
exists a unique up to a constant $U_q \mathfrak{s}(\mathfrak{u}_n \times
\mathfrak{u}_n)$-invariant scalar product. Fix such scalar products via the
integral over the Shilov boundary of the quantum matrix ball
\cite{Cauhy_Sege}:
\begin{equation*}
<u,v>_{\overline{\mathbf k}}\,=\int \limits_{S(\mathbb U)_q}v^{*}u\, d\mu
\qquad u,v \in V_{\overline{\mathbf k}}.
\end{equation*}
Hence each invariant scalar product $(\cdot,\cdot):V \times V \rightarrow
\mathbb C$ is assigned to a set $\{c_{\overline{\mathbf
k}}\}_{\overline{\mathbf k} \in K^{+}} \subset \mathbb{R}_+$ such that
\hbox{$(u,v)=c_{\overline{\mathbf k}}<u,v>_{\overline{\mathbf k}}$} for all
$u,v \in V_{\overline{\mathbf k}}$. Conversely, each
$\{c_{\overline{\mathbf k}}\}_{\overline{\mathbf k} \in K^{+}} \subset
\mathbb{R}_+$ defines a unique sesquilinear Hermitian-symmetric positive
definite $U_q \mathfrak{s}(\mathfrak{u}_n \times \mathfrak{u}_n)$-invariant
form in $V$.

Let us find explicit conditions for the coefficients
$\{c_{\overline{\mathbf k}}\}$ to define the
$U_q\mathfrak{su}_{n,n}$-invariant form.
\par Using the decomposition $U_q\mathfrak{sl}_{2n} \simeq
U_q\mathfrak{p}^{-} \otimes U_q\mathfrak{k} \otimes U_q\mathfrak{p}^{+}$
from Section \ref{tech}and the definitions of $U_q\mathfrak{p}^{+}$ and
$U_q\mathfrak{p}^{-}$, we see that it is sufficient to investigate
invariance of $(\cdot,\cdot)$ under the subspaces $\mathfrak{p}_q^{+}$ and
$\mathfrak{p}_q^{-}$. Moreover, it is enough to prove
$\mathfrak{p}_q^{+}$-invariance of $(\cdot,\cdot)$. We can see that if
$(\cdot,\cdot)$ is $\mathfrak{p}_q^{+}$-invariant, then it is
$\mathfrak{p}_q^{-}$-invariant. Indeed, for all $\eta \in
\mathfrak{p}_q^{-},\, u,v \in V$ we have $(\eta u,v)=\overline{(v,\eta u)}
=\overline{(\eta^{*}v,u)}=(u,\eta^{*}v)$.

Investigate the $\mathfrak{p}_q^{+}$-invariance of the form
$(\cdot,\cdot)$. From results of Section \ref{tech} it follows that
$\pi_{\alpha,\beta} (\mathfrak{p}_q^{+}) (V_{\overline{\mathbf k}}) \subset
\bigoplus \limits_{j=1}^n V_{\bf \overline {k} +e_j}$ . Since the isotypic
components $V_{\overline {\mathbf k}}$ are pairwise orthogonal, one need to
check the invariance in ''non-zero cases'' only (that means for $u \in
V_{\overline {\mathbf k}},\,v \in V_{\bf \overline {k} +
e_j}\,,j=1,\ldots,n$). In this case the invariant conditions are the
following: for all $\xi \in U_q\mathfrak{sl}_{2n},\,u \in V_{\overline
{\mathbf k}},\,v \in V_{\bf \overline {k} + e_j}\,,j=1,\ldots,n,$
\begin{align*}
(P_{\bf \overline {k}+e_j}(\pi_{\alpha,\beta}(\xi)u),v)|_{V_{\bf \overline
{k}+e_j}}=(u,P_{\overline {\mathbf k}}(\pi_{\alpha,\beta}(\xi^{*})
v))|_{V_{\overline {\mathbf k}}},
\end{align*}
where $P_{\overline {\mathbf k}}:V \longrightarrow V_{\overline {\mathbf
k}}$ is an orthogonal projection onto $V_{\overline {\mathbf k}}$. In other
words,
\begin{equation*}
c_{\bf \overline {k}+e_j}<P_{\bf \overline { k}+e_j}
(\pi_{\alpha,\beta}(\xi)u),v>_{\bf \overline {k}+e_j}= c_{\overline
{\mathbf k}}<u,P_{\overline {\mathbf k}}(\pi_{\alpha,\beta}(\xi^{*})
v)>_{\overline {\mathbf k}}.
\end{equation*}

First consider the case $\alpha,\beta \not \in \mathbb Z$. Recall that from
Propositions \ref{a},\,\ref{b},\,\ref{c},\,and \ref{d} it follows that in
$(\mathfrak{p}_q^{-} \oplus \mathfrak{p}_q^{+})\otimes V_{\overline
{\mathbf k}}$ there exist $U_q\mathfrak{k}_{ss}$-highest vectors
$\psi_{j,l}^{\pm},\,j,l=1,\ldots,n$ with weights
$(k_1-k_2,\ldots,k_{j-1}-(k_j \pm e_j), (k_j \pm e_j)
-k_{j+1},\ldots,k_{n-1}-k_n,2k_n+\alpha-\beta,k_{n-1}-k_n,\ldots,(k_{n-l+1}
\mp e_l) - k_{n-l+2}, k_{n-l}-(k_{n-l+1} \mp e_l),\ldots,k_1-k_2)$,
respectively. Define $U_q\mathfrak{k}$-invariant maps
$$
T_{\overline {\mathbf k},j}^{\pm}:(\mathfrak{p}_q^{-} \oplus
\mathfrak{p}_q^{+})\otimes V_{\overline {\mathbf k}} \longrightarrow V_{\bf
\overline {k}\pm e_j}
$$
by their values on the $U_q\mathfrak{k_{ss}}$-highest vectors as follows:
\begin{align*}
 T_{\overline {\mathbf k},j}^{+} (\psi_{j,l}^{+}) =
\begin{cases}
\omega_j(\overline {\mathbf k},q) \cdot v_{\bf \overline {k}+e_j}^{h}
&l=n-j+1;
\\0 &l\ne n-j+1;
\end{cases}
\\T_{\overline {\mathbf k},j}^{-} (\psi_{j,l}^{-}) =
\begin{cases}
\varpi_j(\overline {\mathbf k},q) \cdot v_{\bf \overline {k}-e_j}^{h},
&l=n-j+1;
\\0 &l \ne n-j+1.
\end{cases}
\end{align*} Here $v_{\overline {\mathbf k}}^{h}$, $\varpi_j
(\overline {\mathbf k},q)$ and $\omega_j(\overline {\mathbf k},q)$ are
introduced in Propositions \ref{K_reduce}, \ref{pr_2},\,and \ref{pr_3}.
\begin{lemma}
For all $\xi \in \mathfrak{p}_q^{-} \oplus \mathfrak{p}_q^{+},\,u \in
V_{\overline {\mathbf k}},\,j=1,\ldots,n$ the following holds:
\\\,$P_{\bf \overline {k}+e_j}(\pi_{\alpha,\beta}(\xi)u)=q^{-\beta-n/2}[\beta-k_j+j-1]_q
T_{\overline {\mathbf k},j}^{+} (\xi \otimes u)$;
\\\,$P_{\bf \overline {k}-e_j}
(\pi_{\alpha,\beta}(\xi)u)=q^{\alpha+n/2}[\alpha+k_j+n-j]_q T_{\overline
{\mathbf k},j}^{-} (\xi \otimes u)$.
\end{lemma}
{\bf Proof.} The proof completely repeats the proof of Lemma 9.10 of the
paper \cite{Lee}. \hfill $\square$
\medskip

Using the last lemma, we can rewrite the
$U_q\mathfrak{su}_{n,n}$-invariance condition of the scalar product as
follows: for all $\xi \in \mathfrak{p}_q^{-} \oplus \mathfrak{p}_q^{+},\, u
\in V_{\overline {\mathbf k}},\, v \in V_{\bf \overline
{k}+e_j},j=1,\ldots,n$
\begin{align*}
q^{-\beta-n/2}[\beta-k_j+j-1]_q \,c_{\bf \overline {k}+e_j} <T_{\overline
{\mathbf k},j}^{+}(\xi \otimes u),v>_{\bf \overline {k}+e_j} =
\\=\overline{q^{\alpha+n/2}[\alpha+(k_j+1)+n-j]_q} \,c_{\overline {\mathbf k}}<u,
T_{{\bf \overline {k}+e_j},j}^{-}(\xi^{*} \otimes v)>_{\overline {\mathbf
k}}.
\end{align*}
\begin{proposition}
$<T_{\overline {\mathbf k},j}^{+}(\xi \otimes u),v>_{\bf \overline {k}+e_j}
= -<u,T_{{\bf \overline {k}+e_j},j}^{-}(\xi^{*} \otimes v)>_{\overline
{\mathbf k}}$ for all $j=1,\ldots,n.$
\end{proposition}

{\bf Proof.} Since the maps $T_{\overline {\mathbf k},j}^{\pm}$ does not
depend on $\alpha,\,\beta \in \mathcal{D}$, it is enough to consider only
the special case $\operatorname {Re}(\alpha+\beta)=-n$. In this case the
representation $\pi_{\alpha,\beta}$ is unitary, thus we can put
$c_{\overline {\mathbf k}}=1$ for all $\overline{\mathbf k} \in
\widehat{K}$. Since $\overline{q^{\alpha+n}}=q^{-\beta}$, we see that
\begin{equation*}
q^{-\beta-n/2}\frac{q^{\beta-k_j+j-1}-q^{-\beta+k_j-j+1}}{q-q^{-1}}
(<T_{\overline {\mathbf k},j}^{+}(\xi \otimes u),v>_{\bf \overline {k}+e_j}
+ <u,T_{{\bf \overline {k}+e_j},j}^{-}(\xi^{*} \otimes v)>_{\overline
{\mathbf k}})=0.
\end{equation*}
If we consider non-integral $\alpha,\,\beta$, then $q^{\beta-k_j+j-1}
-q^{-\beta+k_j+1-j}$ does not vanish. This completes the proof.\hfill
$\square$

\medskip

Recall that $\alpha,\beta \not \in \mathbb Z$. Thus the
$U_q\mathfrak{su}_{n,n}$-invariance condition of the scalar product can be
rewritten as follows: for all $\overline{\mathbf k} \in \widehat{K},
j=1,\ldots,n$
\begin{equation}\label{eq_unit}
(1-q^{2(-\beta+k_j+1-j)}) (\overline{1-q^{2(\alpha+(k_j+1)+n-j)}})^{-1}=
-\frac{c_{\overline {\mathbf k}}}{c_{\bf \overline {k}+e_j}}.
\end{equation}
Since the scalar product must be positive definite, we have the following
necessary conditions for the unitarizability of modules of the degenerate
principal series (recall that $q=e^{-h/2}$): for all $\overline{\mathbf k}
\in \widehat{K}, j=1,\ldots,n$
\begin{equation*}
\operatorname{sh} \frac{h}{2}(\beta-k_j+j-1)(\overline{\operatorname{sh}
\frac{h}{2} (\alpha+(k_j+1)+n-j)})^{-1}<0.
\end{equation*}

Using these inequalities, we can present the following series of simple
unitary {\it representations of degenerate principal series} related to the
Shilov boundary.

\medskip

{\it The principal unitary series:} $\mathrm{Re} (\alpha+\beta)=-n$,
$\alpha,\beta \not \in \mathbb Z$. In this case all representations are
unitary. The invariant scalar product provided by the $U_q
\mathfrak{su}_{n,n}$-invariant integral \cite{Cauhy_Sege}.

{\it The complementary series:} $\mathrm{Im} (\alpha+\beta)=0$,
$|\mathrm{Re} \alpha+n|<1$, $|\mathrm{Re} \beta|<1$, $(\mathrm{Re}
\alpha+n)\mathrm{Re} \beta<0$, $\alpha,\beta \not \in \mathbb Z$. In this
case the representations $\pi_{\alpha,\beta}$ are unitary too. (The
required invariant scalar product $(\cdot,\cdot)$ is defined by the
coefficients $\{c_{\overline{\mathbf{k}}}\}$ as follows: let $c_{\bf
\overline{0}}=1$, other coefficients are computed from recurrent relations
such as (\ref{eq_unit}).)

{\it The strange series:} $\mathrm{Im} \alpha=\frac{\pi}{h}$. For such
values of the parameters the respective representations
$\pi_{\alpha,\beta}$ are irreducible and unitary. This series of
representations has no classical analogue. For the first time it appears in
unpublished works of L.Korogodsky and in A.Klimyk and V.Groza's paper (see
\cite{KL_Pak}).

\medskip

Now let $\alpha,\beta \in \mathbb{Z}$. (Recall that in this case
$\pi_{\alpha,\beta}$ is reducible.) For such $\alpha,\beta$ there might
exist unitarizable simple submodules in the respective module (we will
mention them below), although the module is not unitarizable. For each
simple submodule the same arguments as in ''general case'' on the $U_q
\mathfrak{su}_{n,n}$-invariance of scalar product can be applied. In each
case we have the necessary conditions like (\ref{eq_unit}), however they
must be satisfied only on a certain part of $\widehat{K}$. Consider all
possible cases:

{\bf Case 1.} $\alpha+\beta\geq 2-n$. In this case the representation is
not unitary and its unique irreducible subrepresentation is not unitary
too.

{\bf Case 2.} $\alpha+\beta=1-n$. In this case there exist $n$ irreducible
unitary subrepresentations of the representation $\pi_{\alpha,1-n-\alpha}$.
Precisely, $V^s_j$ (see (\ref{c_2_submod})) is a simple submodule in $V$
for any $j=1,\ldots,n$, Notice that each $V^s_j$ can be equipped with a
$U_q \mathfrak{su}_{n,n}$-invariant scalar product $(\cdot,\cdot)$. Such
modules are called small representations because they have ''poor''
decompositions into isotypic components.

{\bf Case 3.} $\alpha+\beta=-n$. In this case the representations are
completely reducible, their irreducible subrepresentations
$V^s_i,\,i=1,\ldots,n+1$ (see Section \ref{reduce}) are unitary (actually,
the required invariant scalar product is the same as for the principal
unitary series).

{\bf Case 4.} $\alpha+\beta \leq -1-n$. In this case the submodules
$V^s_i,\,i=1,\ldots,n+1$ are unitary although there exist non-unitarizable
quotient modules in $V$.

\section{Acknowledgment}

The author is grateful to professor L.L.Vaksman for constant attention to
her work and helpful discussions.

\section {Appendix}

Let us prove Lemma \ref{l_1}. This proof is a q-analogue of the proof of
Lemma 3.4 from \cite{Lee}.

{\bf Proof of Lemma \ref{l_1}.} Statements (\ref{l_1_1})-(\ref{l_1_4}) can
be easily checked.

For example, check the equality $K_jF_{mj}=qF_{mj}K_j$. For $j-m=1$, i.e.,
$m=j-1$, we see that
$K_jF_{j-1,j}=K_jF_{j-1}K_{j-1}=F_{j-1}K_{j-1}K_j=qF_{j-1,j}K_j$. Assume
that for $j-m < r$ equations (\ref{l_1_1})-(\ref{l_1_4}) are proved. Let
$j-m=r$. Then,
\begin{multline*}
\!\!\!K_jF_{mj}=K_j(F_{m+1,j}F_mK_m +\!\!\!\sum_{s=m+2}^{j}
\!\!(-1)^{s+m+1}F_{sj} \ad_{F_{s-1}}\ldots \ad_{F_{m+1}}(F_mK_m)
K(j,m+1,s-1))
\\=K_jF_{m+1,j}F_mK_m+\sum_{s=m+2}^{j} (-1)^{s+m+1}K_j F_{sj} \ad_{F_{s-1}}\ldots
\ad_{F_{m+1}}(F_mK_m) K(j,m+1,s-1)
\\ =qF_{m+1,j}K_jF_mK_m+q\!\!\sum_{s=m+2}^{j-1} \!\!(-1)^{s+m+1}F_{sj} K_j
\ad_{F_{s-1}}\ldots \ad_{F_{m+1}}(F_mK_m) K(j,m+1,s-1)
\\+q(-1)^{j+m+1}\ad_{F_{j-1}}\ldots \ad_{F_{m+1}}(F_mK_m)K_j=qF_{mj}K_j.
\end{multline*}
The proof is completed by induction.

Using (\ref{l_1_4}), prove equality (\ref{l_1_6}). Recall that
$[x]_q=\frac{q^x-q^{-x}}{q-q^{-1}}$. In the next equality we assume that
$Z=E_m$: $E_mF_{mj}=$
\begin{multline*}
= \sum \limits_{s=m+2}^j \!\!\!(-1)^{s+m+1} \!E_m F_{sj}
\ad_{F_{s-1}}\ldots \ad_{F_{m+1}}(F_m K_m) K(j,m+1,s-1) + E_mF_{m+1,j}
F_mK_m
\end{multline*}
\begin{multline*}
\equiv \sum_{s=m+2}^{j} \!\!(-1)^{s+m+1} F_{sj} E_m \ad_{F_{s-1}}\ldots
\ad_{F_{m+1}}(F_mK_m) K(j,m+1,s-1) + q F_{m+1,j} E_m F_m K_m
\\ \equiv q\sum_{s=m+2}^{j} (-1)^{s+m} F_{sj} \ad_{F_{s-1}}\ldots
\ad_{F_{m+2}}(F_{m+1}K_{m+1}) K_m^2 K(j,m+1,s-1) + qF_{m+1,j} [H_m]_q K_m
\\ =q \!\!(\sum_{s=m+2}^{j} \!\!\!(-1)^{s+m+2} F_{sj} \ad_{F_{s-1}}\ldots
\ad_{F_{m+2}}(F_{m+1}K_{m+1}) K(j,m+2,s-1))K_m^2
\\ \cdot q^{j-m-1}K_{m+1}\ldots K_{j-1}[H_{m+1}+\ldots+H_{j-1}+j-m-1]_q) +
qF_{m+1,j} [H_m]_q K_m
\\ = q F_{m+1,j} (q^{j-m-1}K_m^2K_{m+1}\ldots K_{j-1}
[H_{m+1}+\ldots+H_{j-1}+j-m-1]_q +[H_m]_qK_m)
\\ = F_{m+1,j} q^{j-m}K_m\ldots K _{j-1}
[H_m+\ldots+H_{j-1}+j-m-1]_q \; (\!\!\!\!\!\mod U_q \mathfrak{sl}_n \cdot
E_m).
\end{multline*}
Prove equality (\ref{l_1_5}) by induction. If $j-m=2$ and $m<i<j$, then
$i=j-1$, and (\ref{l_1_5}) means that $E_{j-1}F_{j-2,j} \equiv 0 \;
(\!\!\!\!\mod U_q \mathfrak{sl}_n \cdot E_{j-1}).$ It can be proved as
follows:
\begin{multline*}
\!\!\!\!E_{j-1}F_{j-2,j}=E_{j-1}( F_{j-1,j} \,F_{j-2} K_{j-2}
\,K(j,j-1,j-2)- \ad_{F_{j-1}}(F_{j-2} K_{j-2}) \,K(j,j-1,j-1))
\\ \equiv [H_{j-1}]_qK_{j-1} F_{j-2}K_{j-2} - qF_{j-2}K_{j-2} K_{j-1}
[H_{j-1}+1]_q =0 \;(\!\!\!\!\!\mod\, U_q \mathfrak{sl}_n \cdot E_{j-1}).
\end{multline*}
For the inductive step it is sufficient to check that for all $m < i < j$
\begin{multline*}
\!\!\!E_iF_{mj}= \!\!\!\sum \limits_{s=m+2}^j (-1)^{s+m+1} E_i F_{sj}
\ad_{F_{s-1}}\ldots \ad_{F_{m+1}}(F_m K_m) K(j,m+1,s-1) + E_iF_{m+1,j} F_m
K_m
\\ \equiv \sum \limits_{s=m+2}^{i-1} (-1)^{s+m+1} E_i F_{sj}
\ad_{F_{s-1}}\ldots \ad_{F_{m+1}}(F_m K_m) K(j,m+1,s-1) +
E_iF_{m+1,j}F_mK_m
\\ + \sum \limits_{s=i+1}^j \!\!(-1)^{s+m+1} E_i F_{sj}
\ad_{F_{s-1}}\ldots \ad_{F_{m+1}}(F_m K_m) K(j,m+1,s-1)
\\ +(-1)^{i+m+1} E_i F_{ij} \ad_{F_{i-1}}\ldots \ad_{F_{m+1}}(F_m K_m)
K(j,m+1,i-1),
\end{multline*}
(we use (\ref{l_1_4}) and (\ref{l_1_6}) and assume $Z=E_i$). By the
inductive hypothesis, for $s < i$ we have $E_iF_{sj} \equiv 0
\;(\!\!\!\!\mod\, U_q \mathfrak{sl}_n \cdot E_i)$, therefore for all $m < s
< i$ there exists an element $X_s \in U_q \mathfrak{sl}_n$ such that
$E_iF_{sj}=X_sE_i$. From (\ref{l_1_4}), $E_iF_{sj}=F_{sj}E_i$. From
(\ref{l_1_6}), $E_iF_{ij} = qF_{i+1,j}q^{j-i-1}K_i\ldots K_{j-1}
[H_i+\ldots+H_{j-1}+j-i-1]_q.$ Thus,
\begin{multline*}
E_iF_{mj} \equiv \sum \limits_{s=m+2}^{i-1} (-1)^{s+m+1} X_s E_i
\ad_{F_{s-1}}\ldots \ad_{F_{m+1}}(F_m K_m) K(j,m+1,s-1)
\\ +X_{m+1}E_iF_mK_m+ \sum \limits_{s=i+1}^j (-1)^{s+m+1} F_{sj} E_i
\ad_{F_{s-1}}\ldots \ad_{F_{m+1}}(F_m K_m) K(j,m+1,s-1)
\\ +(-1)^{i+m+1} F_{i+1,j} q^{j-i}K_i\ldots K_{j-1}
[H_i+\ldots+H_{j-1}+j-i-1]_q \ad_{F_{i-1}}\ldots \ad_{F_{m+1}}(F_m K_m)
\\ \cdot K(j,m+1,i-1) \equiv (-1)^{i+m+1} q F_{i+1,j}
\ad_{F_{i-1}}\ldots \ad_{F_{m+1}}(F_m K_m) K(j,m+1,i)
\\+(-1)^{i+m} F_{i+1,j} E_i \ad_{F_i}\ldots \ad_{F_{m+1}}(F_m K_m) K(j,m+1,i) =0
\;(\!\!\!\!\!\mod\, U_q \mathfrak{sl}_n \cdot E_i). \tag*{$\square$}
\end{multline*}

The proof of Lemma \ref{l_2} is similar.

Let us prove Proposition \ref{pr_2}. We just have to compute the
coefficients $c_j(\beta,k_j)$. Recall that there is a $U_{q}\mathfrak{sl}_n
\otimes U_{q}\mathfrak{sl}_n$-isomorphism $j_1: \mathfrak{p}_q^{+} \simeq
\mathbb C^n \otimes (\mathbb C^n)^{*},$ where $\mathbb C^n$ is the vector
representation of $U_q \mathfrak{sl}_n$. The isomorphism $j_1^{-1}$ on the
elements of the standard basis for $\mathbb C^n \otimes (\mathbb C^n)^{*}$
is defined as follows:
\begin{eqnarray*}
 &j_1^{-1}\begin{pmatrix} \varepsilon_1 \otimes \varepsilon_1^{*} & \ldots &
\varepsilon_1 \otimes \varepsilon_n^* \\ \ldots & \ldots & \ldots \\
\varepsilon_{n-1} \otimes \varepsilon_1^* & \ldots &\ldots\\ \varepsilon_n
\otimes \varepsilon_1^* & \ldots & \varepsilon_n \otimes \varepsilon_n^*
\end{pmatrix}=
\\ \!\!=&\!\!\!\!\!\!\begin{pmatrix} \ad_{E_1} \ldots \ad_{E_{n-1}}E_n & \ldots & \ldots &
(-1)^{n-1}\ad_{E_{2n-1}}\ldots \ad_{E_{n+1}}\ad_{E_1}\ldots
\ad_{E_{n-1}}E_n
\\ \ldots & \ldots & \ldots & \ldots \\
\ad_{E_{n-1}}E_n & \ldots & \ldots & \ldots \\ E_n & -\ad_{E_{n+1}}E_n &
\ldots & (-1)^{n-1} \ad_{E_{2n-1}} \ldots \ad_{E_{n+1}}E_n
\end{pmatrix}
\end{eqnarray*}
(This follows from the equalities $\ad_{F_j}E_n=0$, $\ad_{E_j}^{2}E_n=0$
for $j=1,\ldots,2n-1,j \neq n$, $\ad_{K_j}E_n=E_n$ for
$j=1,\ldots,n-2,n+2,\ldots,2n-1$, $\ad_{K_j}E_n=q^{-1}E_n$ for $j=n-1$ or
$j=n+1$.) Consider the following embeddings of vector spaces
\begin{align*}
 &\iota_1: U_q \mathfrak{sl}_n \hookrightarrow U_q \mathfrak{sl}_n \otimes
U_q \mathfrak{sl}_n, \qquad \xi \mapsto \xi \otimes 1;
\\& \iota_2: U_q \mathfrak{sl}_n \hookrightarrow U_q \mathfrak{sl}_n \otimes
U_q \mathfrak{sl}_n, \qquad \xi \mapsto 1 \otimes \xi.
\end{align*}
Set $\xi^{(1)}=\iota_1(\xi)$ and $\xi^{(2)}=\iota_2(\xi)$.

From Propositions \ref{a} and \ref{b}, we deduce that for all
$j=1,\ldots,n, \overline{\bf k} \in \widehat{K}$
\begin{multline*}
\mathcal{M^+}_{\overline{\mathbf k}}(\zeta_j \otimes \xi_{n-j+1}) =
\mathcal{M^+}_{\overline{\mathbf k}}(\sum \limits_{m=1}^j (-q^2)^{m-1}
\varepsilon_m \otimes F^{(1)}_{mj}K^{(1)}_-(j,1,m-1)u
\\ \otimes \sum \limits_{m=n-j+1}^n \varepsilon_m^* \, \otimes
 S^{(2)}_{n-j+1,m}\,L^{(2)}_-(n-j+1,m+1,n)\,u^*)
= \mathcal{M^+} _{\overline{\mathbf k}} (\sum \limits_{m=1}^j \sum
\limits_{l=n-j+1}^n (-q^2)^{m-1} \varepsilon_m \\\otimes \varepsilon^*_l
\otimes F_{mj}K_-(j,1,m-1)u \otimes S_{n-j+1,l}\,L_-(n-j+1,l+1,n)\,u^*).
\end{multline*}
\begin{proposition} \label{reverse}
For all $j=1,\ldots,n,\overline{\bf k} \in \widehat{K}$
$$\mathcal{M^+}_{\overline{\mathbf k}} (\zeta_j \otimes \xi_{n-j+1}) =
\lambda^-(n-j+1,n-j+2,n)\mathcal{M^+}_{\overline{\mathbf k}}(\zeta_j
\otimes \varepsilon_{n-j+1}^* \otimes u^*),$$ where
$L^{(2)}_-(n-j+1,n-j+2,n) (v^h_{\overline{\mathbf k}})
=\lambda^-(n-j+1,n-j+2,n) v^h_{\overline{\mathbf k}}$.
\end{proposition}

{\bf Proof.} In the same way as in \cite{Lee}, we have
\begin{multline*}
\mathcal{M^+}_{\overline{\mathbf k}}(\zeta_j \otimes \xi_{n-j+1})
=\mathcal{M^+}_{\overline{\mathbf k}}(\sum \limits_{m=1}^j \sum
\limits_{l=n-j+1}^n (-q^2)^{m-1} \varepsilon_m \otimes \varepsilon^*_l
\otimes F^{(1)}_{mj} K^{(1)}_-(j,1,m-1)u \\\otimes
S^{(2)}_{n-j+1,l}\,L^{(2)}_-(n-j+1,l+1,n)u^*) =
\mathcal{M^+}_{\overline{\mathbf k}} (\sum \limits_{m=1}^j(-q^2)^{m-1}
\varepsilon_m \otimes \varepsilon^*_{n-j+1} \otimes
F^{(1)}_{mj}K^{(1)}_-(j,1,m-1)u \\\otimes L^{(2)}_-(n-j+1,n-j+2+1,n)u^*)
+\mathcal{M^+}_{\overline{\mathbf k}}(\sum \limits_{m=1}^j \sum
\limits_{l=n-j+2}^n (-q^2)^{m-1} \varepsilon_m \\\otimes \varepsilon^*_l
\otimes F^{(1)}_{mj}K^{(1)}_-(j,1,m-1)u \otimes
S^{(2)}_{n-j+1,l}\,L^{(2)}_-(n-j+1,l+1,n)u^*)
\end{multline*}
\begin{multline*}
= \lambda^-(n-j+1,n-j+2,n) \mathcal{M^+}_{\overline{\mathbf k}}(\sum
\limits_{m=1}^j(-q^2)^{m-1} \varepsilon_m \otimes \varepsilon^*_{n-j+1}
\otimes F^{(1)}_{mj}K^{(1)}_-(j,1,m-1)u \otimes u^*)
\\ + \mathcal{M^+}_{\overline{\mathbf k}} (\sum \limits_{m=1}^j \sum
\limits_{l=n-j+2}^n \!\!\!(-q^2)^{m-1} \varepsilon_m \otimes
\varepsilon^*_l \otimes F^{(1)}_{mj} K^{(1)}_-(j,1,m-1)u \otimes
S^{(2)}_{n-j+1,l} L^{(2)}_-(n-j+1,l+1,n)u^*)
\\ =\lambda^-(n-j+1,n-j+2,n)
\mathcal{M^+}_{\overline{\mathbf k}} (\zeta_j \otimes \varepsilon_{n-j+1}^*
\otimes u^*)
\\ +\mathcal{M^+}_{\overline{\mathbf k}} (\sum \limits_{l=n-j+2}^n \sum
\limits_{m=1}^j (-q^2)^{m-1} \varepsilon_m \otimes (\ad_{E_{n-l-2}}\ldots
\ad_{E_j})^{(2)}\varepsilon^*_{n-j+1}
\\ \otimes F^{(1)}_{mj}
K^{(1)}_-(j,1,m-1)u \otimes S^{(2)}_{n-j+1,l} L^{(2)}_-(n-j+1,l+1,n)\,u^*)
\\ =\lambda^-(n-j+1,n-j+2,n)
\mathcal{M^+}_{\overline{\mathbf k}} (\zeta_j \otimes \varepsilon_{n-j+1}^*
\otimes u^*)
\\ +\sum_{l=n-j+1}^n (\ad_{E_{n-l-2}}\ldots \ad_{E_j})^{(2)} S^{(2)}_{n-j+1,l}
L^{(2)}_-(n-j+1,l+1,n) \mathcal{M^+}_{\overline{\mathbf k}} (\zeta_j
\otimes \varepsilon_{n-j+1}^* \otimes u^*).
\end{multline*}
The vector $\mathcal{M^+}_{\overline{\mathbf k}} (\zeta_j \otimes
\varepsilon_{n-j+1}^* \otimes u^*) \in V_{\bf \overline{k}+e_j}$ and is a
$U_q \mathfrak{sl}_n \otimes 1$-highest vector. Therefore
$\mathcal{M^+}_{\overline{\mathbf k}} (\zeta_j \otimes
\varepsilon_{n-j+1}^* \otimes u^*) \in V_{\bf \overline{k}+e_j}$,
$\mathcal{M^+}_{\bf \overline{k}} (\zeta_j \otimes \varepsilon_{n-j+1}^*
\otimes u^*)=c \cdot v^h_{\bf \overline{k}+e_j}$ with some $c \in \mathbb
C$. Now we conclude that in the obtained expression all summands except the
first equal 0. \hfill $\square$

To find $c_j(\beta,k_j)$ we must compute
\begin{multline}\tag*{(A1)}\label{re}
\mathcal{M^+}_{\overline{\mathbf k}} (\zeta_j \otimes \varepsilon_{n-j+1}^*
\otimes u^*)= \mathcal{M^+}_{\overline{\mathbf k}}(\sum \limits_{m=1}^j
(-q^2)^{m-1} \varepsilon_m \otimes \varepsilon^*_{n-j+1} \otimes
F^{(1)}_{mj} K^{(1)}_-(j,1,m-1)u \otimes u^*)
\\ =\sum \limits_{m=1}^j (-q^2)^{m-1}
\pi_{\alpha,\beta} \!(\!(-1)^{j-1} \ad_{E_{n+j-1}}\ldots \ad_{E_{n+1}}
\ad_{E_m}\ldots \!\ad_{E_{n-1}}E_n)
\\F^{(1)}_{mj}K^{(1)}_-(j,1,m-1)(v^h_{\overline{\mathbf k}}).
\end{multline}
We need some auxiliary lemmas. Recall that in this paper we introduce the
notation for q-minors of the matrix $\mathbf z$ (see (\ref{q_minz})). Set
$\mathbf z^{\wedge k}_{a_1,\ldots,a_k}= \mathbf z^{\wedge k
\{1,\ldots,k\}}_{\,\,\,\,\,\,\{a_1,\ldots,a_k\}}.$
\begin{lemma}\label{l_min}
For all $1 \leq m \leq k \leq j-2$
\begin{multline*}
(-q)^{j-k-1}\mathbf z^{\wedge j-1} \mathbf z^{\wedge
k}_{1,\ldots,m-1,m+1,\ldots,j} - \sum \limits_{s=k+1}^{j-2}(-q)^{s-k-1}
\mathbf z^{\wedge j-1}_{1,\ldots,s-1,s+1,\ldots,j} \mathbf z^{\wedge
k}_{1,\ldots,m-1,m+1,\ldots,s}
\\=\mathbf z^{\wedge
j-1}_{1,\ldots,m-1,m+1,\ldots,j} \mathbf z^{\wedge k}. \hfill \square
\end{multline*}
\end{lemma}
\begin{lemma} \label{G}
For all $1 \leq m \leq j \leq n$ we have $F_{mj}=q^{j-m-1}G_{mj}$, where
\begin{equation*}
G_{mj}=F_mK_mF_{m+1,j}+\!\!\!\!\sum_{s=m+2}^{j} \!\!\!(-q)^{s-m-1}
\ad_{F_{s-1}}\ldots \ad_{F_{m+1}}(F_mK_m) F_{sj} K_-(j,m+1,s-1). \hfill
\square
\end{equation*}
\end{lemma}
\begin{lemma} \label{FG}
For all $1 \leq m \leq j \leq n$ $$G_{mj}(v^h_{\overline{\mathbf k}}) =
(q^{1/2})^{j-m} \kappa_-(j,m,j-1) z^{\wedge
j-1}_{1,\ldots,m-1,m+1,\ldots,j} \frac{v^h_{\overline{\mathbf
k}}}{z^{\wedge j-1}},$$ where $K_-(j,m+1,j-1)(v^h_{\overline{\mathbf k}})=
\kappa_-(j,m,j-1) v^h_{\overline{\mathbf k}}.$
\end{lemma}

{\bf Proof.} We prove this lemma by induction. For $j-m=1$ the statement is
obvious, since
\begin{multline*}
G_{m,m+1}(v^h_{\overline{\mathbf k}})=F_m K_m ((z^{\wedge 1})^{k_1-k_2}
\ldots(z^{\wedge n})^{k_n})
\\ = q^{k_m-k_{m+1}} q^{1/2} [k_m-k_{m+1}]_q (z^{\wedge
1})^{k_1-k_2}\ldots(z^{\wedge m-1})^{k_{m-1}-k_m} z^{\wedge
m}_{1,\ldots,m-1,m+1} (z^{\wedge m})^{k_m-k_{m+1}-1}
\\ \cdot
(z^{\wedge m+1})^{k_{m+1}-k_{m+2}}\ldots (z^{\wedge n})^{k_n} = q^{1/2}
\kappa_-(m+1,m,m) z^{\wedge m}_{1,\ldots,m-1,m+1} v^h_{\bf
\overline{k}-e_m}.
\end{multline*}
For the proof of the inductive step we use two previous lemmas. By Lemma
\ref{G}, we have
\begin{multline*}
G_{mj}(v^h_{\overline{\mathbf k}}) = \sum_{s=m+2}^{j} (-q)^{s-m-1}
\ad_{F_{s-1}}\ldots \ad_{F_{m+1}}(F_mK_m) \cdot F_{sj} K_-(j,m+1,s-1))
(v^h_{\overline{\mathbf k}})
\\ + F_mK_mF_{m+1,j} (v^h_{\overline{\mathbf
k}})= F_mK_mF_{m+1,j} (v^h_{\overline{\mathbf k}})
\\ + \sum_{s=m+2}^{j} (-q)^{s-m-1} \kappa_-(j,m+1,s-1)
\ad_{F_{s-1}}\ldots \ad_{F_{m+1}}(F_mK_m) F_{sj} (v^h_{\overline{\mathbf
k}}).
\end{multline*}
By the inductive hypothesis, for all $j-s<j-m$
\begin{multline*}
G_{mj}(v^h_{\overline{\mathbf k}}) = q^{1/2(j-m-1)} \kappa_-(j,m+1,j-1)
F_mK_m ( z^{\wedge j-1} _{1,\ldots,m-1,m+1,\ldots,j}
\frac{v^h_{\overline{\mathbf k}}}{z^{\wedge j-1}}) + \sum \limits_{s=m+2}^j
(-q)^{s-m-1}
\\ \cdot q^{1/2(j-s)} \kappa_-(j,m+1,s-1) \kappa_-(j,s,j-1)
\ad_{F_{s-1}}\ldots \ad_{F_{m+1}}(F_mK_m) (z^{\wedge j-1}
_{1,\ldots,m-1,m+1,\ldots,j} \frac{v^h_{\overline{\mathbf k}}}{z^{\wedge
j-1}})
\\ =q^{1/2(j-m-1)}\kappa_-(j,m+1,j-1) \cdot ( F_mK_m
(z^{\wedge j-1}_{1,\ldots,m-1,m+1,\ldots,j} \frac{v^h_{\overline{\mathbf
k}}} {z^{\wedge j-1}})
\\ +\sum \limits_{s=m+2}^j (-q^{1/2})^{s-m-1} \ad_{F_{s-1}}\ldots \ad_{F_{m+1}}
(F_mK_m) (z^{\wedge j-1}_{1,\ldots,s-1,s+1,\ldots,j}
\frac{v^h_{\overline{\mathbf k}}} {z^{\wedge j-1}}).
\end{multline*}
Using the explicit formulas for the $U_q\mathfrak{sl}_{2n}$-action in
$\mathbb C [\mathrm{Mat}_n]_q$ and properties of the comultiplication (see
Section \ref{probl}), we obtain that
\begin{multline*}
\ad_{F_{s-1}}\ldots \ad_{F_{m+1}} (F_mK_m) (z^{\wedge
j-1}_{1,\ldots,s-1,s+1,\ldots,j} \frac{v^h_{\overline{\mathbf k}}}
{z^{\wedge j-1}})
\\ =(q^{1/2})^{s-m} (q^{k_m-k_{m+1}}[k_m-k_{m+1}]_q
(z^{\wedge 1})^{k_1-k_2} \cdot\ldots \cdot (z^{\wedge m-1})^{k_{m-1}-k_m}
z^{\wedge m}_{1,\ldots,m-1,m+1,\ldots,s}
\\ \cdot (z^{\wedge m})^{k_m-k_{m+1}-1}\ldots(z^{\wedge s-1})^{k_{s-1}-k_s}
z^{\wedge j-1}_{1,\ldots,s-1,s+1,\ldots,j} (z^{\wedge s})^{k_s-k_{s+1}}
\cdot\ldots\cdot(z^{\wedge n})^{k_n}
\\ +(-q)q^{k_m-k_{m+1}+k_{m+1}-k_{m+2}}[k_{m+1}-k_{m+2}]_q (z^{\wedge
1})^{k_1-k_2} \cdot\ldots \cdot (z^{\wedge m})^{k_m-k_{m+1}}
\\ \cdot z^{\wedge m+1}_{1,\ldots,m-1,m+1,\ldots,s} (z^{\wedge
m+1})^{k_{m+1}-k_{m+2}-1}\ldots(z^{\wedge s-1})^{k_{s-1}-k_s} z^{\wedge
j-1}_{1,\ldots,s-1,s+1,\ldots,j}
\\ \cdot (z^{\wedge s})^{k_s-k_{s+1}} \cdot\ldots\cdot (z^{\wedge n})^{k_n}+ \ldots+
(-q)^{s-2}q^{k_m-k_{m+1}+\ldots+k_{s-2}-k_{s-1}}[k_{s-2}-k_{s-1}]_q
\\ \cdot (z^{\wedge 1})^{k_1-k_2} \cdot\ldots \cdot (z^{\wedge s-2})^{k_{s-2}-k_{s-1}}
z^{\wedge s-1}_{1,\ldots,m-1,m+1,\ldots,s} \cdot (z^{\wedge
s-1})^{k_{s-1}-k_s-1} z^{\wedge j-1}_{1,\ldots,s-1,s+1,\ldots,j}
\\ \cdot (z^{\wedge s})^{k_s-k_{s+1}}\cdot\ldots\cdot (z^{\wedge n})^{k_n}+
(-q)^{s-1} q^{k_m-k_{m+1}+\ldots+k_{s-1}-k_s}[k_{s-1}-k_s]_q (z^{\wedge
1})^{k_1-k_2} \cdot\ldots
\\ \cdot (z^{\wedge s-2})^{k_{s-2}-k_{s-1}}
(z^{\wedge s-1})^{k_{s-1}-k_s} z^{\wedge j-1}_{1,\ldots,m-1,m+1,\ldots,j}
(z^{\wedge s})^{k_s-k_{s+1}}\cdot\ldots\cdot(z^{\wedge n})^{k_n}).
\end{multline*}
Finally, by Lemma \ref{l_min}, we see that
\begin{multline*}
G_{mj}(v^h_{\overline{\mathbf k}}) =(q^{1/2})^{j-m}\kappa_-(j,m+1,j-1)
\cdot z^{\wedge j-1}_{1,\ldots,m-1,m+1,\ldots,j}
\frac{v^h_{\overline{\mathbf k}} }{z^{\wedge j-1}}\cdot
\frac{q^{2(j-2+k_1-k_j)}-1}{q-q^{-1}}. \tag*{$\square$}
\end{multline*}
By the last lemma and (\ref{re}), we can compute the coefficients
$c_j(\beta,k_j)$ introduced in Proposition \ref{pr_2}.
\begin{proposition}\label{21}
For all $1 \leq j \leq n$ $$\mathcal{M^+}_{\overline{\mathbf k}} (\zeta_j
\otimes \varepsilon_{n-j+1}^* \otimes u^*)=
q^{-\beta-n/2+k_j+j}[\beta-k_j+j-1]_q \kappa_-(j,1,j-1)v^h_{\bf
\overline{k}+e_j}.$$
\end{proposition}
{\bf Proof} We have
\begin{multline*}
\!\!\!\!\mathcal{M^+}_{\overline{\mathbf k}} (\zeta_j \otimes
\varepsilon_{n-j+1}^* \otimes u^*)= \sum \limits_{m=1}^j (-q^2)^{m-1}
\pi_{\alpha,\beta} ((-1)^{j-1}\ad_{E_{n+j-1}}\ldots \ad_{E_{n+1}}
\ad_{E_m}\ldots \ad_{E_{n-1}}E_n)
\\ \cdot F^{(1)}_{mj}K^{(1)}_-(j,1,m-1) v^h_{\overline{\mathbf k}}=
\sum \limits_{m=1}^j (-q^2)^{m-1} \kappa_-(j,1,m-1)
\\ \cdot
\pi_{\alpha,\beta} ((-1)^{j-1}\ad_{E_{n+j-1}}\ldots \ad_{E_{n+1}}
\ad_{E_m}\ldots \ad_{E_{n-1}}E_n) F^{(1)}_{mj} v^h_{\overline{\mathbf k}}.
\end{multline*}
By Lemma \ref{FG},
\begin{multline*}
\mathcal{M^+}_{\overline{\mathbf k}} (\zeta_j \otimes \varepsilon_{n-j+1}^*
\otimes u^*)= \sum \limits_{m=1}^{j-1} (-q^2)^{m-1} q^{j-m-1}
(q^{1/2})^{j-m} \kappa_-(j,m,j-1) \kappa_-(j,1,m-1)
\\ \cdot (-1)^{j-1}\pi_{\alpha,\beta}
(\ad_{E_{n+j-1}}\ldots \ad_{E_{n+1}} \ad_{E_m}\ldots \ad_{E_{n-1}}E_n)
(z^{\wedge j-1}_{1,\ldots,m-1,m+1,\ldots,j} \frac{v^h_{\overline{\mathbf
k}}}{z^{\wedge j-1}})
\\ +(-q^2)^{j-1} \kappa_-(j,1,j-1) \pi_{\alpha,\beta} ((-1)^{j-1}
\ad_{E_{n+j-1}}\ldots \ad_{E_{n+1}} \ad_{E_j}\ldots \ad_{E_{n-1}}E_n)
(v^h_{\overline{\mathbf k}})
\\ =  q^{3/2j-3} \kappa_-(j,1,j-1) \sum \limits_{m=1}^{j-1} (-1)^{j+m} q^{m/2}
\\ \cdot \pi_{\alpha,\beta} (\ad_{E_{n+j-1}}\ldots \ad_{E_{n+1}}
\ad_{E_m}\ldots \ad_{E_{n-1}}E_n) (z^{\wedge
j-1}_{1,\ldots,m-1,m+1,\ldots,j} \frac{v^h_{\overline{\mathbf
k}}}{z^{\wedge j-1}})
\\ + q^{2j-2} \kappa_-(j,1,j-1) \pi_{\alpha,\beta} (\ad_{E_{n+j-1}}\ldots
\ad_{E_{n+1}} \ad_{E_j}\ldots \ad_{E_{n-1}}E_n) (v^h_{\overline{\mathbf
k}}).
\end{multline*}

In Section \ref{probl} the following morphism of
$U_q\mathfrak{sl}_{2n}$-modules was defined:
$$
\iota: \mathbb C[\mathrm{Mat}_n]_q \rightarrow \mathbb{C}
[\mathrm{Pl}_{n,2n}]_{q,t}, \qquad \iota(z^{\wedge
k\{b_1,\ldots,b_k\}}_{\,\,\,\,\{a_1,\ldots,a_k\}}) =t^{-1}t^{\wedge
n}_{\{1,\ldots,n\}J},
$$
with $J=\{n+1,\ldots,2n\} \setminus \{ 2n+1-b_1,\ldots,2n+1-b_k\} \cup
\{a_1,\ldots,a_k\}$. Therefore $\iota(z^{\wedge k})= t^{-1} t^{\wedge
n}_{\{1,\ldots,n\}\{1,\ldots,k,n+1,\ldots,2n-k\}}$. It follows that
$$
\iota(v^h_{\overline{\mathbf k}})=q^{c}t^{-k_1}(t^{\wedge
n}_{\{1,\ldots,n\}\{1,n+1,\ldots,2n-1\}})^{k_1-k_2}\ldots(t^{\wedge
n}_{\{1,\ldots,n\}\{1,\ldots,n\}})^{k_n},
$$
where some $c \in \mathbb C$. Using the definition of $\pi_{\alpha,\beta}$,
we obtain that for all $\xi \in U_q\mathfrak{sl}_{2n}$, $f \in \mathbb{C}
[\mathrm{Pl}_{n,2n}]_{q,t}$
\begin{equation*}
\pi_{\alpha,\beta} (\xi) (f)= q^{c} t^{-\beta} \xi \cdot (t^{\beta} f
(t^{\wedge n}_{\{1,\ldots,n\}\{1,\ldots,n\}})^{\alpha})(t^{\wedge
n}_{\{1,\ldots,n\}\{1,\ldots,n\}})^{-\alpha}.
\end{equation*}
For $m<j$
\begin{multline*}
(\ad_{E_{n+j-1}}\ldots \ad_{E_{n+1}} \ad_{E_m}\ldots \ad_{E_{n-1}}E_n)
(t^{\wedge n
\{1,\ldots,n\}}_{\,\,\,\,\{1,\ldots,m-1,m+1,\ldots,j,n+1,\ldots,2n-j,2n-m+1\}})
\\ = (-1)^{m-j}(q^{-1/2})^{j+n-m} q^{m-j-1} t^{\wedge n}
_{\{1,\ldots,n\}\{1,\ldots,j,n+1,\ldots,2n-j\}},
\end{multline*}
and for $m=j$ $(\ad_{E_{n+j-1}}\ldots \ad_{E_{n+1}} \ad_{E_j}\ldots
\ad_{E_{n-1}}E_n) (t^{\wedge
n}_{\{1,\ldots,n\}\{1,\ldots,j-1,n+1,\ldots,2n-j+1\}})= (q^{-1/2})^n
t^{\wedge n}_{\{1,\ldots,n\}\{1,\ldots,j,n+1,\ldots,2n-j\}}$. For other
summands we use an analogue of Lemma \ref{l_min}. Finally, we have
\begin{multline*}
\mathcal{M^+}_{\overline{\mathbf k}} (\zeta_j \otimes \varepsilon_{n-j+1}^*
\otimes u^*)= \kappa_-(j,1,j-1)q^{-n/2}\cdot (q^2\sum \limits_{m=1}^{j-1}
q^{2m-2}+ \frac{1-q^{-2\beta+2k_1}}{1-q^{-2}}
\\+ \sum_{m=1}^{j-1}
q^{-2\beta+2k_m}\frac{1-q^{-2k_m+2k_{m+1}}}{1-q^{-2}})\cdot v^h_{\bf
\overline{k}+e_j}= q^{-\beta-n/2+k_j+j} \kappa_-(j,1,j-1) [\beta-k_j+j-1]_q
v^h_{\bf \overline{k}+e_j}. \tag*{$\square$}
\end{multline*}

Repeat the same arguments to prove Proposition \ref{pr_3}. First, we have
the explicit formulas for the isomorphism $j_2: \mathfrak{p}_q^{-} \simeq
(\mathbb C^n)^* \otimes \mathbb C^n$:
\begin{eqnarray*}
&j_2^{-1} \begin{pmatrix} \varepsilon_1^* \otimes \varepsilon_1 & \ldots &
\varepsilon_n^* \otimes \varepsilon_1 \\ \ldots & \ldots & \ldots \\
\varepsilon_1^* \otimes \varepsilon_n & \ldots &\varepsilon_n^* \otimes
\varepsilon_n
\end{pmatrix}=
\\ &\!\!\!\!\!\!\!\!\!\begin{pmatrix} \!\!\!\!\!\!(-1)^{n-1} \ad_{F_1} \ldots \ad_{F_{n-1}}(K_nF_n)
& \!\!\!\ldots & \!\!\!\!-\ad_{F_{n-1}}(K_nF_n) &K_nF_n
\\ \ldots & \!\!\!\ldots & \!\!\!\!\ldots & \ad_{F_{n+1}}(K_nF_n) \\ \ldots \\
(-1)^{n-1} \ad_{F_{2n-1}}\ldots \ad_{F_{n+1}}\ad_{F_1}\ldots
\ad_{F_{n-1}}(K_nF_n) & \!\!\!\ldots & \!\!\!\!\ldots &
\!\!\!\!\!\!\!\!\!\!\!\!\!\!\!\!\!\!\! \ad_{F_{2n-1}}\ldots
\ad_{F_{n+1}}(K_nF_n)
\end{pmatrix}
\end{eqnarray*}
For the proof of Proposition \ref{pr_3} we must compute the following:
\begin{multline*}
\mathcal{M^-}_{\overline{\mathbf k}}(\xi'_j \otimes \zeta'_{n-j+1}) =
\mathcal{M^-}_{\overline{\mathbf k}}(\sum \limits_{m=j}^n \varepsilon_m^*
\otimes S^{(1)}_{jm}L^{(1)}_-(j,m+1,n)u \otimes \sum \limits_{m=1}^{n-j+1}
(-q^2)^{m-1} \varepsilon_m
\\ \otimes F^{(2)}_{m,n-j+1}K^{(2)}_-(n-j+1,1,m-1)u^*) =
\mathcal{M^-}_{\overline{\mathbf k}} (\sum \limits_{m=j}^n \sum
\limits_{l=1}^{n-j+1} (-q^2)^{l-1} \varepsilon_m^* \otimes \varepsilon_l
\\ \otimes S^{(1)}_{jm}L^{(1)}_-(j,m+1,n)u \otimes F^{(2)}_{l,n-j+1}
K^{(2)}_-(n-j+1,1,l-1)u^*).
\end{multline*}
\begin{proposition}
For all $1 \leq j \leq n$ $$\mathcal{M^-}_{\overline{\mathbf k}}(\xi'_j
\otimes \zeta'_{n-j+1}) =
\lambda_-(j,j+1,n)\mathcal{M^-}_{\overline{\mathbf k}}(\varepsilon^*_j
\otimes u \otimes \zeta'_{n-j+1}),$$ where
$L^{(1)}_-(j,j+1,n)u=\lambda_-(j,j+1,n)u$.
\end{proposition}
The proof is similar to the proof of Proposition \ref{reverse}.

Therefore in order to find the coefficients $d_j(\alpha,k_j)$ introduced in
Proposition \ref{pr_3} we must only compute
\begin{multline*}
\mathcal{M^-}_{\overline{\mathbf k}}(\varepsilon^*_j \otimes u \otimes
\zeta'_{n-j+1})
\\=\mathcal{M^-}_{\overline{\mathbf k}}(\sum \limits_{l=1}^{n-j+1}
(-q^2)^{l-1} \varepsilon_j^* \otimes \varepsilon_l \otimes u \otimes
F_{l,n-j+1}K_-(n-j+1,1,l-1)u^*)
\\ =\sum \limits_{l=1}^{n-j+1} (-q^2)^{l-1}
\pi_{\alpha,\beta}(\ad_{F_j}\ldots \ad_{F_{n-1}}\ad_{F_{n-1+l}}\ldots
\ad_{F_{n+1}}(K_nF_n))
\\ \cdot F^{(2)}_{l,n-j+1}K^{(2)}_-(n-j+1,1,l-1)
(v^h_{\overline{\mathbf k}}).
\end{multline*}
These computations are analogous to the ones from the proof of Proposition
\ref{21}.

\newpage
\onecolumn

\newpage
\begin{picture} (300,310)(0,0)
\put(15,20){\dottedline{2}(0,0)(0,300)}
\put(15,325){\vector(0,1){5}} \put(5,30){\dottedline{2}(0,0)(300,0)}
\put(310,30){\vector(1,0){5}} \put(5,100){\line(1,0){310}}
\put(5,180){\line(1,0){310}} \put(65,20){\line(0,1){310}}
\put(125,20){\line(0,1){310}} \put(65,40){\vector(1,0){15}}
\put(300,180){\vector(0,-1){15}} \put(300,100){\vector(0,1){15}}
\put(125,40){\vector(-1,0){15}} \put(10,25){\line(1,1){295}}
\put(5,305){$k_2$} \put(305,20){$k_1$} \put(70,20){$-\alpha-1$}
\put(130,20){$\beta$} \put(0,105){$-\alpha$} \put(0,185){$\beta+1$}

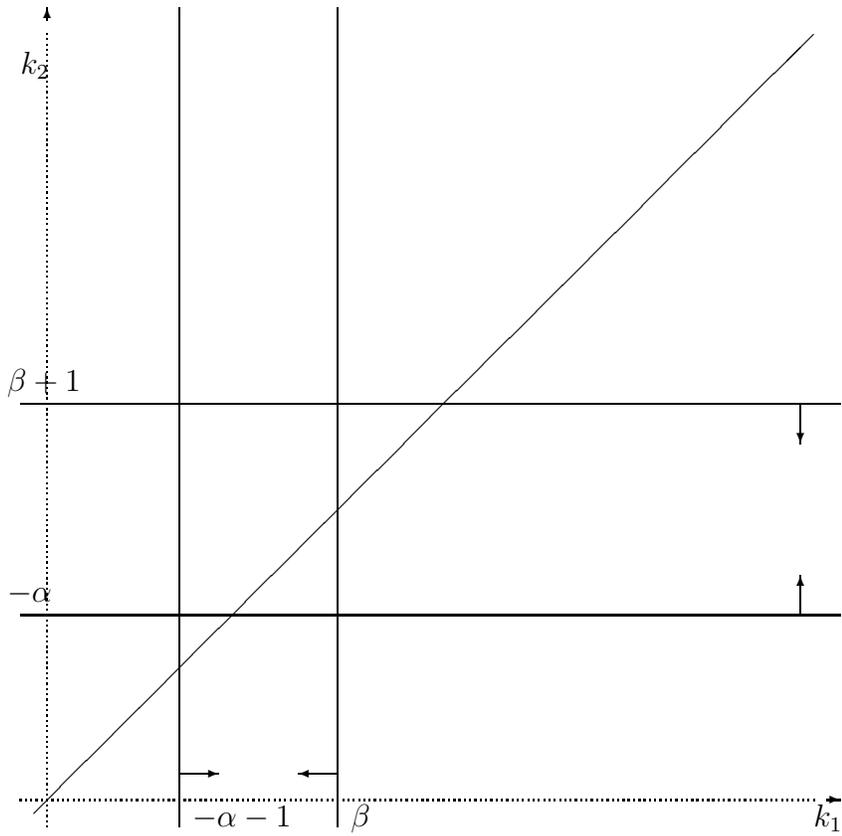
\put(100,0){Fig.1. Structure of $\pi_{\alpha,\beta}$ with
$\alpha+\beta \geq 0$.}
\end{picture}

\begin{picture} (300,360)(0,0)
\put(15,20){\dottedline{2}(0,0)(0,300)}
\put(15,315){\vector(0,1){5}} \put(5,30){\dottedline{2}(0,0)(300,0)}
\put(310,30){\vector(1,0){5}} \put(5,130){\line(1,0){310}}
\put(105,10){\line(0,1){310}} \put(90,40){\vector(1,0){15}}
\put(300,145){\vector(0,-1){15}} \put(300,115){\vector(0,1){15}}
\put(120,40){\vector(-1,0){15}} \put(10,25){\line(1,1){295}}
\put(5,310){$k_2$} \put(305,20){$k_1$} \put(110,20){$\beta$}
\put(0,135){$-\alpha$}

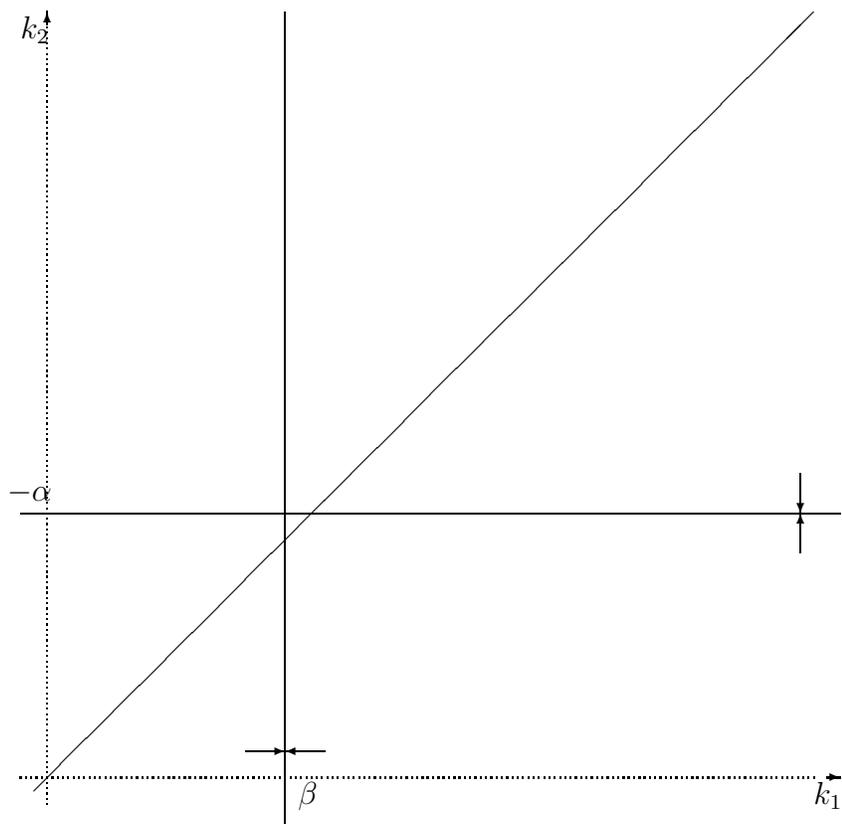
\put(100,0){Fig.2. Structure of $\pi_{\alpha,\beta}$ with
$\alpha+\beta=-1$.}
\end{picture}

\begin{picture}(300,310)(0,0)
\put(15,20){\dottedline{2}(0,0)(0,300)}
\put(15,325){\vector(0,1){5}} \put(5,30){\dottedline{2}(0,0)(300,0)}
\put(310,30){\vector(1,0){5}} \put(5,120){\line(1,0){310}}
\put(5,200){\line(1,0){310}} \put(65,20){\line(0,1){310}}
\put(145,20){\line(0,1){310}} \put(145,40){\vector(1,0){15}}
\put(285,120){\vector(0,-1){15}} \put(285,200){\vector(0,1){15}}
\put(65,40){\vector(-1,0){15}} \put(10,25){\line(1,1){295}}
\put(5,320){$k_2$} \put(305,20){$k_1$} \put(70,20){$\beta$}
\put(150,20){$-\alpha-1$} \put(0,125){$\beta+1$}
\put(0,205){$-\alpha$}

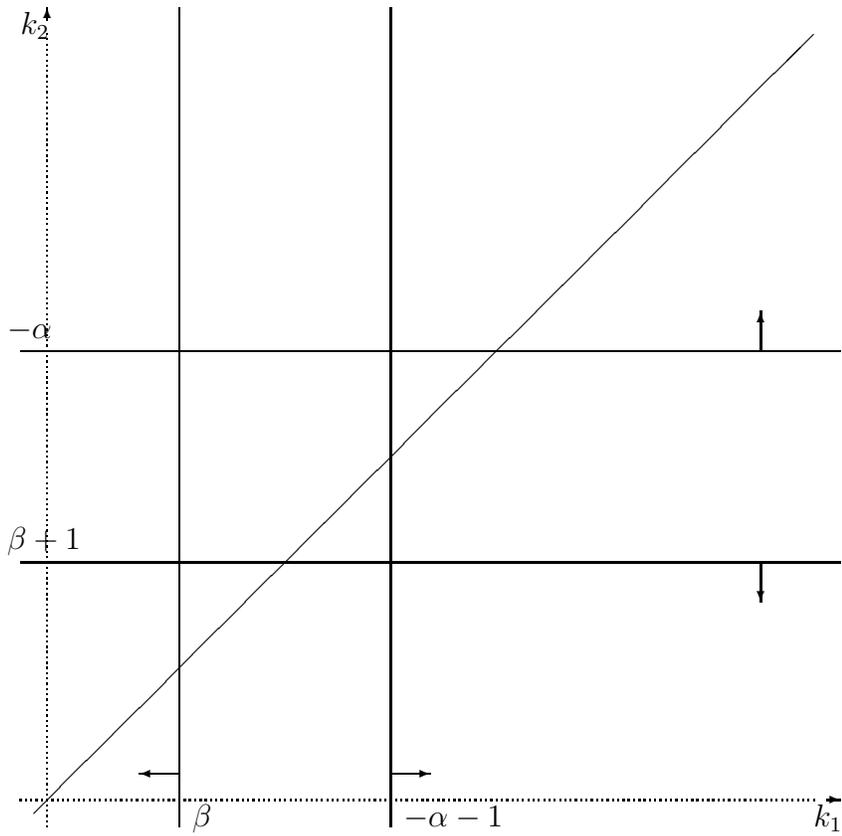
\put(100,0){Fig.3. Structure of $\pi_{\alpha,\beta}$ with
$\alpha+\beta=-2$.}
\end{picture}

\begin{picture} (300,360)(0,0)
\put(15,20){\dottedline{2}(0,0)(0,300)}
\put(15,325){\vector(0,1){5}} \put(5,30){\dottedline{2}(0,0)(300,0)}
\put(305,30){\vector(1,0){5}} \put(5,100){\line(1,0){310}}
\put(5,180){\line(1,0){310}} \put(65,20){\line(0,1){310}}
\put(125,20){\line(0,1){310}} \put(125,40){\vector(1,0){15}}
\put(300,100){\vector(0,-1){15}} \put(300,180){\vector(0,1){15}}
\put(65,40){\vector(-1,0){15}} \put(10,25){\line(1,1){295}}
\put(5,320){$k_2$} \put(305,20){$k_1$} \put(70,20){$\beta$}
\put(130,20){$-\alpha-1$} \put(0,105){$\beta+1$}
\put(0,185){$-\alpha$}

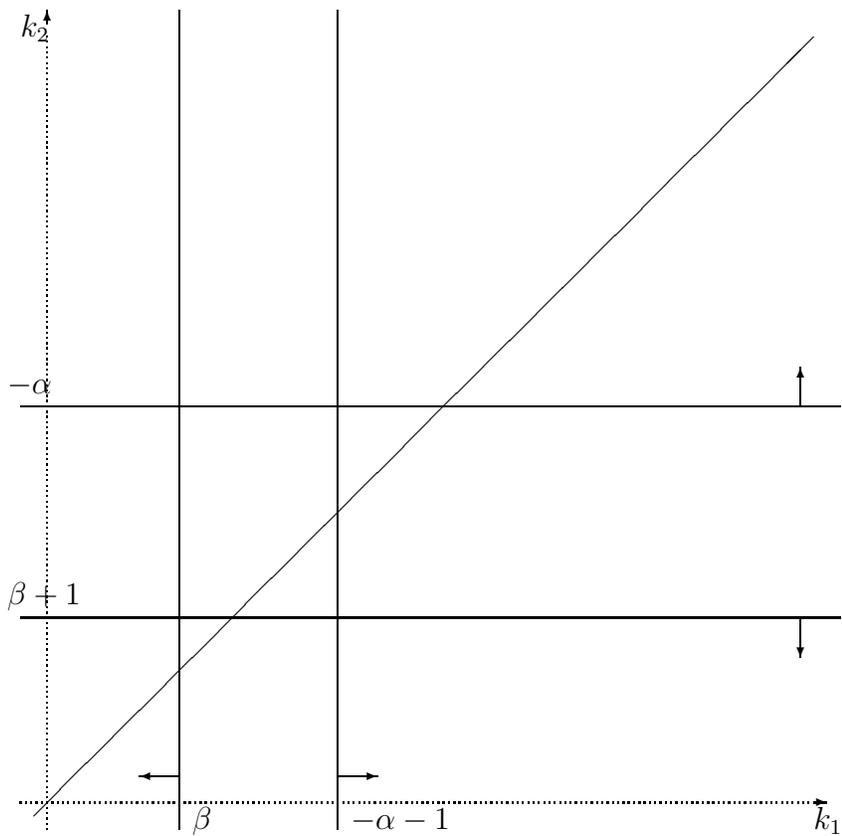
\put(100,0){Fig.4. Structure of $\pi_{\alpha,\beta}$ with
$\alpha+\beta \leq -3$.}
\end{picture}
\end{document}